\documentclass{article}
\usepackage{xcolor}
\usepackage{amsmath,amsfonts,amssymb}
\usepackage{ntheorem}
\usepackage[all]{xy}
\usepackage{mathrsfs}
\usepackage[margin=3cm]{geometry}
\usepackage{graphicx}
\usepackage{paralist}

\numberwithin{equation}{section}

\theoremstyle{plain}
\theoremheaderfont{\normalfont\bfseries}
\theorembodyfont{\slshape}
\theoremseparator{.}
\newtheorem{teo}{Theorem}[section]
\newtheorem{lem}[teo]{Lemma}
\newtheorem{p}[teo]{Proposition}
\newtheorem{cor}[teo]{Corollary}

\theoremstyle{plain}
\theoremheaderfont{\normalfont\bfseries}
\theorembodyfont{\upshape}
\theoremseparator{:}
\newtheorem{obs}[teo]{Remark}

\theoremstyle{plain}
\theorembodyfont{\upshape}
\theoremseparator{:}
\newtheorem{defi}[teo]{Definition}
\newtheorem{exa}[teo]{Example}

\theoremstyle{change}
\theorembodyfont{\upshape}
\theoremseparator{:}
\newtheorem{ej}[teo]{Example}
    \theoremclass{ej}
    \theoremstyle{changebreak}

   \theoremstyle{change}
\theorembodyfont{\upshape}
\theoremseparator{:}
    \theoremclass{ex}
    \theoremstyle{changebreak}

    \theoremheaderfont{\sc}
    \theorembodyfont{\upshape}
    \theoremstyle{nonumberplain}
    \theoremseparator{.}
    \newtheorem{dem}{Proof}

\newcommand{\fin}{\hfill$\blacksquare$}
\newcommand{\mb}{\mathbb}

\newcommand{\stackbounds}[2]{\stackrel{\mbox{\scriptsize \(#1\)}}{\mbox{\scriptsize \(#2\)}}}
\newcommand{\codimC}{\mathop{\mathrm{codim}_\mathbb{C}}}
\newcommand{\comp}{\mathop{\mathrm{comp}}}

\newcommand{\Gr}{\mathbb{G}}
\newcommand{\setpart}[1]{#1_{\rm set}}

\title{Matroids and the space of torus-invariant subvarieties of the Grassmannian
with given homology class}
\author{E. Javier Elizondo\footnote{Thanks to UNAM and the DGAPA for their support through the fellowship PASPA.} 
\and Alex Fink\footnote{Partially funded by the European Union's Horizon 2020 research and innovation programme under the Marie Sk\l odowska-Curie grant agreement No~792432.}
\and Cristhian Garay L\'opez\footnote{Supported by CONACYT through Project 299261 and a postdoctoral scholarship.}}
\date{ \today}

\begin{document}
\maketitle
\begin{abstract}
Let $\Gr(d,n)$ be the complex Grassmannian of affine $d$-planes in $n$-space. 
We study the problem of characterizing the set of algebraic subvarieties of $\Gr(d,n)$
invariant under the action of the maximal torus $T$
and having given homology class $\lambda$. 
We give a complete answer for the case where $\lambda$ is the class of a $T$-orbit, and partial results for other cases,
using techniques inspired by matroid theory.
This problem has applications to the computation of the Euler-Chow series for Grassmannians of projective lines. 
\end{abstract}

\section{Introduction}
\label{Section_Intro}
The problem that inspired this article is the computation of the $p$-dimensional Euler-Chow series for $\Gr(d,n)$, the complex Grassmannian of affine $d$-planes in $n$-space. This series can be defined for any projective variety as follows.

Given a projective variety $X$, denote by $M$ the monoid in $H_{2p}(X;
\mathbb{Z})$ generated by algebraic classes of effective cycles. $M$ is
considered as a multiplicative monoid: if $a, b$ are effective cycles, we
set $t^{a}t^{b}=t^{a+b}$, where “$t$” is a formal variable that turns
addition into multiplication. The set of functions from $M$ to the
integers is denoted by $\mathbb{Z}[\![M]\!]$ and the Euler-Chow series of $X$
is defined as
\begin{equation}
\label{eq_E_p}
    E_p(X):=\sum_{\lambda\in M}\chi(\mathcal{C}_{p, \lambda}(X))t^\lambda \in\mathbb{Z}[\![M]\!],
\end{equation}
where $\mathcal{C}_{p, \lambda}$ is the Chow variety parametrizing effective
algebraic $p$-cycles, and $\chi(\mathcal{C}_{p, \lambda})$ is its topological Euler characteristic.
See Section~\ref{Section_Link} for a fuller introduction.

Computing the Euler characteristic of Chow varieties is a very hard problem. 
Introduced as merely a way to record the computations of individual Euler characteristics, 
it was later seen that the Euler-Chow series has many interesting properties: for instance, in some cases it agrees with the Hilbert series of a variety. A fundamental paper where some formulae for projective bundles are shown and examples were computed for the series together with formal definitions is \cite{eli-lim}. 
For the case of simplicial toric varities the series was computed in~\cite{eli-tor}, and it is also known for the case of abelian varieties \cite{eli&hain-abe}. 
There is also a motivic version of the series and some interesting examples are computed  in \cite{eli-kim}, where a relation is noted between the Weil zeta series and the motivic version of the Euler-Chow  series. In \cite{xi-javier} a relation between the Cox ring and the Euler-Chow series is shown.  

There is a natural action of the  maximal algebraic torus $T \cong (\mathbb{C}^{*})^{n}$ on $ \mb{G}(d,n)$, 
and it is known \cite{law&yau-hosy} that $\chi(\mathcal{C}_{p,\lambda}(\mb{G}(d,n))) \, = \, \chi(\mathcal{C}_{p,\lambda}(\mb{G}(d,n))^{T})$. 
The main results in this paper are
parametrizations of the invariant subvarieties of $ \mb{G}(d,n)$ under the action of~$T$ 
whose homology class is $\lambda$.
This in principle enables us to compute the Euler characteristics of Chow varieties that appear in~\eqref{eq_E_p}.
In some cases we are able to carry out this computation.

Although the general case is very hard, we were able to give an answer when the invariant subvarieties are closures of $T$-orbits, or homologous to these. For the particular case of $ \mb{G}(2,n)$ we are able to go further. 
The main techniques used in this paper come from the theory of matroids, and the fact that any matroid of rank two is representable over $\mathbb{C}$ allowed us to exploit the relation between invariant subvarieties and matroids for $ \mb{G}(2,n)$.

\noindent\textbf{Acknowledgements.} We thank the referee for their close reading, especially for their simplifications of multiple proofs.

\subsection{Statement of results}
In Theorem \ref{teo_orbits} we show that any $T$-invariant irreducible subvariety  of $\Gr(d,n)$ which is homologous to the closure of a $T$-orbit 
is in fact itself a closure of a $T$-orbit. 
As a result, we obtain a concrete description of the space that parametrizes $T$-orbits with a fixed homology class. 
Then we show in Theorem \ref{theo_complete_set} how to construct a complete set of representatives $\Pi_n$ for the orbit space of the action of the symmetric group $\mathfrak{S}_n$ on the set of thin Schubert cells of $\mb{G}(2,n)$.
In Theorem \ref{thm_map} we show that the set of classes of the closure of the thin Schubert cells is precisely the set of all possible monomials in the Schubert classes. 
In Theorem \ref{teo_decompodecompo} an effective algorithm to compute cohomology classes of torus orbits in $\mb{G}(2,n)$.
Last, in Theorem~\ref{C_3d G24}, we compute the Euler-Chow series of $\mb{G}(2,4)$.

\subsection{Organization of the article}
\label{Section_Roadmap}
The first part of this article works on the general case $\Gr(d,n)$.
In the hope that this article will be of interest to algebraic geometers as well as to combinatorialists, we start by introducing some preliminaries from both fields in Section \ref{Section_Preliminaries}, like matroids (Section \ref{Subsection_Matroids}) and the thin Schubert cell decomposition of the Grassmannians (Section \ref{Subsection_Decomposition}). Particularly important for us is the bundle structure that we get on each thin Schubert cell of $\mb{G}(d,n)$, when we take its  quotient by the torus $T$. 

Section \ref{Section_T_Primitive} is the technical part of this work. Here we study the problem of characterizing the set of algebraic subvarieties of 
$\Gr(d,n)$ invariant under the action of $T$ and with a given class $\lambda\in H_{*}(\Gr(d,n), \mathbb{Z})$. 
In Section \ref{Subsection_Classes_Of_Orbits} we give a complete answer when $\lambda$ is the class of (the closure of) a $T$-orbit, namely, 
any invariant subvariety homologous to an orbit must be an orbit.

The second part of the article starts in Section  \ref{Section_case_lines} and focuses on the case $\mb{G}(2,n).$  In Section \ref{Sect_geo} we discuss the geometry of the thin Schubert cells considered as torus bundles. 
Prior to this we introduce the language of set partitions which we use to work with matroids of rank~2.

In Section \ref{Section_Ring_part} we recast some aspects of the ring structure of $H_{*}(\Gr(d,n), \mathbb{Z})$ in terms of partial partitions of the integer $n$. We do this for Schubert classes in Section \ref{Section_Schubert} and the general case in Section \ref{Subsection_Thin_SC}.
This leaves us well positioned to compute Euler-Chow series for Grassmannians $\mb{G}(2,n)$.
As a proof of concept we exhibit the computations for $\mb{G}(2,4)$, the intereseting case being the 3-cycles, in Section~\ref{Section_Link}.
We intend further such computations to be the subject of future work.


\section{Preliminaries}
\label{Section_Preliminaries}
\paragraph{Conventions.} By a variety $X$ we mean an integral algebraic scheme over $\mathbb{C}$, and by a subvariety of $ X$ we mean an integral closed subscheme. Thus  we will work  with irreducible closed algebraic subsets $Y\subset X$, and if some scheme structure has to be invoked, it will be its reduced closed subscheme structure (unless otherwise stated).

Throughout $0<d<n$ will be two integers. We will denote by $[n]$ the set $\{1,\ldots,n\}$ frequently endowed with its usual order,  and by $\binom{[n]}{d}$ the set of subsets of cardinality $d$ drawn from $[n]$, also endowed with its  term-wise partial order  \cite[p. 8]{BGW}:
for $I=(i_1,\ldots,i_d), J=(j_1,\ldots,j_d)\in \binom{[n]}{d}$  such that $i_1<\cdots<i_d$ and $j_1<\cdots<j_d$, 
we have that $I\leq J$ if and only if $i_k\leq j_k$ for all $1\leq k\leq d$.
We denote by  $\mathfrak{S}_n$  the symmetric group on $n$ letters. 

If $S$ is a set, we denote by $S\setminus T$ the complement of $T$ within~$S$. If $S$ is a finite set, we denote its number of elements by $\#S$. Sometimes we write $S\sqcup T$ to indicate disjoint union of sets.


A \textbf{decomposition} of a topological space $X$ is a family $\{S_i\}_{i\in P}$ of pairwise disjoint and locally-closed subsets of $X$, satisfying $X=\bigsqcup_{i\in P} S_i$. We say that the decomposition $\{S_i\}_{i\in P}$ is a \textbf{stratification} if the closure  $\overline{S}_i$ of each piece is itself a union of pieces. In this case the set of pieces becomes a partially ordered set if we define $i\leq j$ if and only if $S_i\subset\overline{S}_j$.

\subsection{Matroids}
\label{Subsection_Matroids} 
Matroids are combinatorial objects suitable for recording the discrete linear-algebraic information carried by some $x\in \Gr(d,n)$,
and they are the indexing object for the thin Schubert decomposition of Proposition~\ref{P_stratification_1} that is central to this work.
This section introduces them.

Matroids are known for admitting many equivalent definitions --
Gian-Carlo Rota called this phenomenon by Garrett Birkhoff's name ``cryptomorphism'' --
and we use multiple of them.
In general part of the data of a matroid is a choice of a set of elements. We will always take this set to be $[n]$.

\begin{defi}
\label{def_matroid_bases}
A \textbf{matroid} (of \textbf{rank} $d$, on $n$ elements) is a pair $M=([n],\mathcal B(M))$, where $\mathcal B(M)\subset \binom{[n]}{d}$ is a non-empty family which satisfies the following exchange property: 
\begin{equation}
\label{mep}
\text{for every }I,J\in \mathcal B(M)\text{ and }i\in I\setminus J,\text{ there exists }j\in J\setminus I\text{ such that }(I\setminus\{i\})\cup\{j\}\in \mathcal B(M).
\end{equation}
\end{defi}
The  elements of $\mathcal B(M)$ are the \textbf{bases} of $M$, and $[n]$ is called the \textbf{ground set}.
We denote by $\mathcal{M}^d_n$ the set of matroids of rank~$d$ on $n$ elements.

We now introduce some standard terms and constructions in matroid theory, for which we  recommend \cite{JO} as a reference. 
Let $M$ be a matroid. An \textbf{independent set} of $M$ is a set $I\subset [n]$ that is contained in some basis of $M$. 
A set which is not independent is said to be \textbf{dependent}. 
A \textbf{circuit} of $M$ is a minimal dependent subset. If $S\subset[n]$, we denote $M|S$ the \textbf{restriction} of $M$ to $S$: this is the matroid over $S$ whose independent sets are the independent sets of $M$ that are contained in $S$. If $T=[n]\setminus S$, we denote $M\setminus S:=M|T$  the \textbf{deletion} of $S$.
%
We say that $i\in [n]$ is a \textbf{loop} (respectively a \textbf{co-loop}) of $M$ if $i\notin I$ for all $I\in \mathcal B(M)$ (respectively $i\in I$ for all $I\in \mathcal B(M)$). We denote $L(M)$ (respectively $C(M)$)  the set of loops (respectively of co-loops)  of $M$.

Every matroid $M$ has a decomposition into \textbf{connected components}. 
We construct a relation $\sim_M$ on $[n]$ as follows: we say that $i\sim_M j$ if and only if there exist $I,J\in \mathcal B(M)$ such that $J=(I\setminus\{i\})\cup\{j\}$. 
Then $\sim_M$ is an equivalence relation
whose equivalence classes we call the \textbf{connected components} of $M$. 
We denote by  $\comp(M)$ the set of connected components of~$M$.

If $\comp(M)=1$ we say that $M$ is \textbf{connected}. Note that singleton loops and singleton co-loops of~$M$ are connected components.

Matroids can also be characterized as certain lattice polytopes.  
Let $\{e_1,\ldots,e_n\}$ be the standard basis of $\mathbb{R}^n$ as a $\mathbb{R}$-vector space. For every $I\in\binom{[n]}{d}$, set $e_I:=\sum_{i\in I}e_i$. Note that each $e_I$ belongs to the affine space $\{(x_1,\ldots,x_n)\in \mathbb{R}^n\::\:\sum_ix_i=d\}$ of dimension $n-1$. 
Set $\Delta(d,n):=\text{Conv}\{e_I\::\:I\in\binom{[n]}{d}\}\subset \mathbb{R}^n$, which also has dimension $n-1$. 
If  $M=([n],\mathcal B(M))$ is a matroid, then its corresponding \textbf{matroid polytope} is $\Delta(M)=\text{Conv}\{e_I\::\:I\in 
\mathcal B(M) \}$. Conversely, if $\Delta(M)$ is a matroid polytope, then $\mathcal B(M)=\{\text{Supp}(e_I)=I\::\:e_I\in \Delta(M)\}$ satisfies condition \eqref{mep}. 

\begin{teo}[\cite{GGMS}]
A lattice polytope $\Delta\subset \Delta(d,n)$ is a matroid polytope  if and only if $M\neq\emptyset$ and every edge of $\Delta$ is in the direction $e_i-e_j$ for some $i,j\in[n]$.
\end{teo}

\begin{obs}
\label{dimension_of_poly}
The dimension of the matroid polytope $\Delta(M)$ is  $n-\#\comp(M)$. For $0\leq p\leq n-1$, we define $\mathcal{M}^d_n(p)=\{M\in \mathcal{M}^d_n\::\:\dim\Delta(M)=p\}$. 
\end{obs}


\subsection{Thin Schubert cells and matroidal decompositions}
\label{Subsection_Decomposition} 

A standard reference for this section is \cite{BGW}.
We denote by $\mathbb{C}^n$ the set of closed points of $\text{Spec}(\mathbb{C}[x_1,\ldots,x_n])$, and by $\Gr(d,n)$ the complex Grassmannian of $d$-dimensional vector subspaces of $\mathbb{C}^n$. This is a smooth projective variety of complex dimension $d(n-d)$,
with the Grassmann-Pl\"ucker projective embedding $\Gr(d,n)\subset \mathbb{P}^{N(d,n)}$, 
where $N(d,n)=\binom{n}{d}-1$ and $\mathbb{P}^{N(d,n)}$ has homogeneous coordinates $\{x_I\::\:I\in \binom{[n]}{d}\}$.

We denote the \textbf{Schubert varieties} within $\Gr(d,n)$ with respect to a complete flag $\mathcal{V}$ as $\Sigma^\mathcal{V}_{a}$.
As indexing objects we use sequences $a=(a_1,\ldots,a_d)\in\mb{Z}^d$ such that
\begin{equation}
    \label{indexing_sequence}
    n-d\geq a_1\geq\cdots\geq a_d\geq0.
\end{equation}
If $a=(a_1,\ldots,a_d)\in\mb{Z}^d$ satisfies \eqref{indexing_sequence}, we denote by $\mathcal{P}(a)\in\binom{[n]}{d}$ the strictly increasing $d$-tuple $\mathcal{P}(a)=(n-d+1-a_1,n-d+2-a_2,\ldots)$.
Concretely, if $\mathcal{V}$ is the standard flag $(\{0\}=V_0\subset V_1\subset\cdots\subset V_n=\mathbb{C}^n)$, where $V_i=V(x_j\::\:j>i)\cong\mathbb{C}^i$,
then $\Sigma^\mathcal{V}_{a}=\Gr(d,n)\cap H_{a}$, where  $H_{a}\subset \mathbb{P}^{N(n,d)}$ is the coordinate subspace defined by the vanishing of the coordinates $\{x_I\::\:I\nleq \mathcal{P}(a)\}$; see  \cite[Theorem 4.3]{EH}. 

We recall that a $\mb{Z}$-basis for $H^*(\Gr(d,n),\mb{Z})$
is given by the \textbf{Schubert classes} $\sigma_a$,
the cohomology classes of the Schubert varieties. 
The multiplicative structure of $H^*(\Gr(d,n),\mb{Z})$
can be described by the ring isomorphism
\begin{equation}\label{symmetric_functions}
H^*(\Gr(d,n),\mb{Z})\stackrel\sim\to
\mb{Z}[x_1,\ldots,x_d]^{\mathfrak{S}_d} / \langle s_a(x_1,\ldots,x_d) : a_1>n-d\rangle,
\end{equation}
where $s_a(x_1,\ldots,x_d)$ is a Schur symmetric function
(for a reference on symmetric functions we recommend \cite[Chapter~7]{RS2}).
The isomorphism is given by $\sigma_a\mapsto s_a(x_1,\ldots,x_d)$.
In Section~\ref{Subsection_Classes_Of_Orbits} we will use the perspective whereby $H^*(\Gr(d,n),\mb{Z})$ is a quotient of what Stanley calls the algebra $\Lambda$ of symmetric functions.
Stanley's $\Lambda$ is spanned by elements $s_a$ that are notionally symmetric functions in countably many variables $x_i$ as opposed to the $d$ variables in \eqref{symmetric_functions}, 
so the indexing object $a$ is now an infinite nonincreasing list of nonnegative integers that eventually becomes~$0$.
This perspective allows us to view $H^*(\Gr(d,n),\mb{Z})$ as a quotient ring of $H^*(\Gr(d',n'),\mb{Z})$, with a canonical splitting as abelian groups whenever $d\leq d'$ and $n-d\leq n'-d'$.

Let $T$ be a maximal torus acting on $\Gr(d,n)$. For any $x\in \Gr(d,n)$, we have that the closure of the $T$-orbit $\overline{T x}$ is a proper and normal \cite[p. 16]{EHL} toric variety  of complex dimension $0\leq p\leq n-1$.

For $M\in \mathcal{M}_n^d$, the \textbf{thin Schubert cell} is the locally closed locus $G_M\subset\Gr(d,n)$ 
where the coordinate $x_I$ is nonzero if and only if $I\in \mathcal B(M)$. 
The ideal sheaf locally generated by the $x_I$ for $I\not\in\mathcal B(M)$ gives a scheme structure to the closure of $G_M$, and thus to $G_M$ itself; this scheme is not necessarily reduced.
The cell $G_M$ can be empty (in which event matroid theorists say that $M$ is not \textbf{realizable}, or \textbf{representable}, over~$\mathbb{C}$).
Since for any $x\in\Gr(d,n)$ the set $\{I:x_I\ne0\}$ is the set of bases of some matroid,
the next proposition follows.
\begin{p}
\label{P_stratification_1}
We have the following decomposition: 

\begin{equation}
	\label{intro_tscdesc}
\Gr(d,n)=\bigsqcup_{M\in \mathcal{M}_n^d}G_M.
\end{equation}
 
\end{p}

\begin{obs}\label{obs_tscdesc}
The decomposition \eqref{intro_tscdesc} fails to provide a stratification of $\Gr(d,n)$ in general since the closure of a thin Schubert cell is not always a union of thin Schubert cells \cite{EK}. 
However, \eqref{intro_tscdesc} does provide a stratification for $\Gr (2,n)$ (see Proposition~\ref{P_stratification}). 
\end{obs}

The \textbf{Schubert matroid} $M(\Sigma^\mathcal{V}_{a})$ is the matroid of the generic point of $\Sigma^\mathcal{V}_{a}$,
so that $\Sigma^\mathcal{V}_{a}$ is the closure of $G_{M(\Sigma^\mathcal{V}_{a})}$.
By the paragraph after \eqref{indexing_sequence}, we have
\begin{equation}\label{eq_general_Schubert_matroid}
\mathcal B(M(\Sigma^\mathcal{V}_{a})) = \{I\in\binom{[n]}d:I\le\mathcal P(a)\}.
\end{equation}



Each $G_M$ is $T$-invariant, so we  have a  fiber bundle 
\begin{equation}
\label{intro_fb}
T_M\xrightarrow{}G_M\xrightarrow{}\mathcal{G}_M
\end{equation}
where  $T_M$  is the quotient of $T$ by its stabilizer on points of $G_M$ (which depends only on~$M$) and $\mathcal{G}_M:=G_M/T_M$. We take the topological quotient in the analytic topology. 
These bundles will be a key tool from Section~\ref{Section_case_lines} onward.

The action of $\mathfrak{S}_n$ on the set $\{1,\ldots,n\}$ induces an action on $\binom{[n]}{d}$,
descending to an action on $\mathcal{M}_n^d$: an element $\sigma\in \mathfrak{S}_n$ sends the matroid $M=([n],\mathcal B(M))$ to the matroid $\sigma(M)=([n],\sigma(\mathcal B(M)))$. 
It also induces an action of $\mathfrak{S}_n$ on~$\mathbb{C}^n$ by permuting the indices of its coordinates $\{x_1,\ldots,x_n\}$, and thus
an action on $d$-dimensional subspaces of $\mathbb{C}^n$, to wit on~$\Gr(d,n)$.
This last action $\mathfrak{S}_n\circlearrowright\Gr(d,n)$ transforms the coordinates 
$\{x_I\::\:I\in \binom{[n]}{d}\}$ as $\sigma(x_I) = \pm x_{\sigma(I)}$.
Thus it permutes the cells of the thin Schubert decomposition \eqref{intro_tscdesc}: 
$$G_M\cong \sigma(G_M)=G_{\sigma(M)} \text{ for all }\sigma\in \mathfrak{S}_n\text{ and } M\in \mathcal{M}_n^d.$$
In Section~\ref{Section_case_lines} we will construct a complete set of representatives $M_n^d\subset \mathcal{M}_n^d$ for the orbit space  $\mathcal{M}_n^d/\mathfrak{S}_n$ and study  the resulting family of bundles $\mathcal{B}=\{T_{M}\xrightarrow{}G_{M}\xrightarrow{}\mathcal{G}_{M}\:|\:M\in M_n^d\}$. 

\subsection{Matroid subdivisions and valuations}

\begin{defi}
Let $\Delta,\Delta_1,\ldots,\Delta_k\in \mathcal{M}^d_n$ be matroid polytopes of the same dimension in~$\mathbb{R}^n$. 
We say that $\{\Delta_1,\ldots,\Delta_k\}$ is a {\bf matroidal subdivision} of $\Delta$ 
if $\bigcup_{i=1}^k\Delta_i=\Delta$ 
and $\Delta_1,\ldots,\Delta_k$ are the maximal cells of a polyhedral complex in which every cell is a matroid polytope.
\end{defi}

\begin{defi}
\label{def_rigid_matroid}
We say that  $M\in \mathcal{M}^d_n$ is {\bf rigid} if  there are no matroid subdivisions of $\Delta_M$.
\end{defi}

Let  $\{\Delta_1,\ldots,\Delta_k\}$ be a matroidal subdivision of $\Delta$. For $\emptyset\neq S\subset [k]$, we denote $\Delta_S=\bigcap_{i\in S}\Delta_i$, and $\Delta_\emptyset:=\Delta$. 
By definition, every $\Delta_S$ is a matroid if it is nonempty.

Suppose that $M$ is a matroid such that $\mathcal{G}_M$ is positive-dimensional, so it contains a curve $\mathcal{C}$.
After a change of base to a large enough field extension, $\mathcal{C}$ will contain infinitely many points.
By a folklore theorem one of us has previously attributed to Lafforgue \cite{BakerLorscheid}, $M$ is not rigid.
This implies that if $M\in \mathcal{M}_n^d$ is rigid, then $\text{dim}_{\mb{C}}G_M\leq n-1$. 
It is not known whether the converses to these results hold.
For a cautionary example see \cite[Example 3.2]{BakerLorscheid},
a failure of the converse to a related statement involving the inner Tutte group.
Also, curves in $\mathcal{G}_M$ are associated to so-called \textbf{regular} matroidal subdivisions in the sense of \cite[\S2.2.3]{dLRS}, but in principle a matroid polytope may admit only non-regular subdivisions.


As an example of rigidity, a matroid $M$ on ground set $[n]$ is {\bf series-parallel} if, up to isomorphism, either it is $U_2^1$,
or it is obtained from a series-parallel matroid $M'$ on ground set $[n-1]$ by either of the two dual operations {\bf parallel extension} or {\bf series extension}. 
Each of these operations takes an element $i\in[n-1]$ as input.
We say $M$ is a parallel extension of~$M'$ if
\[\mathcal B(M)=\mathcal B(M')\cup\{B\setminus\{i\}\cup\{n\} : B\in\mathcal B(M'),i\in B\}\]
and a series extension of~$M'$ if
\[\mathcal B(M)=\{B\cup\{n\}:B\in\mathcal B(M')\}\cup\{B\cup\{i\} : B\in\mathcal B(M'),i\not\in B\}.\]

\begin{p}\label{p:sp rigid}
Series-parallel matroids are rigid.
\end{p}

\begin{dem}
The proof is inductive on the series-parallel construction process.  As a base case, $U_2^1$ is rigid.
Viewing elements of $\Gr(d,n)$ as representing lists of $n$ vectors spanning $\mb{C}^d$,
parallel extension adjoins a new vector which must be parallel to a given one of the existing vectors:
there is only a $\mb{C}^*$ worth of choices for the new vector, so rigidity is preserved.
Series extension is the dual operation under the duality $\Gr(d,n)\cong \Gr(n-d,n)$, so it also preserves rigidity.\fin
\end{dem}

\begin{defi}
A map $f$ from $\mathcal{M}^d_n$ to an abelian group $G$ is a \textbf{matroid valuation} if 
\begin{equation}\label{matroid_valuation}
\sum_{S\subset [k]}(-1)^{\#S}f(\Delta_S)=0
\end{equation}
whenever $\{\Delta_1,\ldots,\Delta_k\}$ is a matroidal subdivision of $\Delta\in\mathcal{M}^d_n$.

A map $f$ from $\mathcal{M}^d_n(p)$ to $G$ is an \textbf{additive matroid valuation} if
the extension of $f$ to $\bigcup_{q\leq p}\mathcal{M}^d_n(q)$ given by $f(M)=0$ for $\dim\Delta(M)<p$
satisfies \eqref{matroid_valuation} whenever $\{\Delta_1,\ldots,\Delta_k\}$ is a matroidal subdivision of $\Delta\in\mathcal{M}^d_n(p)$.
\end{defi}
\begin{obs}
The additive group generated by indicator functions of matroid polytopes has a filtration by the dimension of the polytope,
with successive subquotients being free abelian groups \cite[Section 8]{DF}.   
By choosing a splitting of this filtration,
any additive matroid valuation in our sense may in fact be extended to a matroid valuation defined on all of $\mathcal{M}^d_n$.
\end{obs}

\section{On the homology of $T$-equivariant subvarieties}
\label{Section_T_Primitive}

If $Y\subset \Gr(d,n)$ is a subvariety, we denote by $h(Y)\in H_*(\Gr(d,n),\mathbb{Z})$ its homology class.
We denote by $\Gr(d,n)^T$ the set  consisting of the subvarieties of $\Gr(d,n)$ which are $T$-invariant. We will be interested in the map 

\begin{equation}
    \label{Homology_map}
    h:\Gr(d,n)^T\xrightarrow{}H_*(\Gr(d,n),\mathbb{Z})
\end{equation}


The problem of computing Euler--Chow series motivates us to look for
\begin{enumerate}
    \item the image of $h$, which we call the set of prime $T$-classes, and
    \item for every prime $T$-class  $\lambda$, the set  $h^{-1}(\lambda)\subset \Gr(d,n)^T$ of $T$-invariant subvarieties of $\Gr(d,n)$ which are homologous to $\lambda$.
\end{enumerate}

The approach we take to answering the latter question is encapsulated in Corollary \ref{cor_param_space} below.
Recall that $\mathcal{G}_M:=G_M/T_M$ for $M\in\mathcal{M}_n^d$, as in~\eqref{intro_fb}. 

\begin{defi}
\label{def_map_V}
If $Y\subset \mathcal{G}_M$ is a subvariety, we can form the pullback bundle $Y\times_{\mathcal{G}_M}G_M$ and define $\mathbb{V}_M(Y):=\overline{Y\times_{\mathcal{G}_M}G_M}$.
\end{defi}


The following result uses the operations $\mathbb{V}_M$ to describe the set $\Gr(d,n)^T$.

\begin{teo}\label{V(Y)} 
We have  $$\Gr(d,n)^T=\bigsqcup_{M\in \mathcal{M}_n^d} \{\mathbb{V}_M(Y) : Y \text{ is a subvariety of } \mathcal{G}_M \}. $$
\end{teo}
\begin{dem}
For each $M\in \mathcal{M}_n^d$, 
$\mathbb{V}_M(\mathcal{G}_M)$ is included in~$\Gr(d,n)^T$ as a subvariety of~$G_M$.
If $Y\subset \Gr(d,n)$ is a  subvariety, 
then there exists a unique $N\in \mathcal{M}_n^d$
and a subvariety $Y_N\subset G_N$ such that $Y=\overline{Y_N}$. 
This $N$ is the matroid whose bases are the indices of Pl\"ucker coordinates which are nonzero at the generic point of~$Y$.
If $Y$ is $T$-invariant, then it induces a subvariety  $Y_N/T=\mathcal{Y}\subset \mathcal{G}_N$, and it follows that $Y=\mathbb{V}_N(\mathcal{Y})$.\fin

\end{dem}

\begin{teo}\label{V(Y)_2} 
Let $M_n^d\subset \mathcal{M}_n^d$ be a complete family of representatives for $\mathcal{M}_n^d/\mathfrak{S}_n$. Then the map \eqref{Homology_map}  factors through the projection  $$\bigcup_{M\in \mathcal{M}_n^d}\{\mathbb{V}_M(Y) : Y \text{ is a subvariety of } \mathcal{G}_M \}\longrightarrow \bigcup_{M\in {M}_n^d}\{\mathbb{V}_M(Y) : Y \text{ is a subvariety of } \mathcal{G}_M \}.$$
\end{teo}

\begin{dem}
The symmetric group $\mathfrak{S}_n$ on $[n]=\{1,\ldots,n\}$ acts by isomorphisms on $\Gr(d,n)$. Letting it act trivially on $H_*(\Gr(d,n),\mathbb{Z})$ makes the map $h$ equivariant: 
that is, if $\sigma\in\mathfrak{S}_n$ and $Y\in \Gr(d,n)^T$, then $h(\sigma(Y))=h(Y)=\sigma(h(Y))$.

In fact $\sigma(G_M)=G_{\sigma(M)}$, so if $M_n^d\subset \mathcal{M}_n^d$ is  a complete family of representatives  for $\mathcal{M}_n^d/\mathfrak{S}_n$, then it follows by Theorem \ref{V(Y)} that, as sets,
the orbit space $\Gr(d,n)^T/\mathfrak{S}_n$ equals 
$\bigcup_{M\in {M}_n^d}\{\mathbb{V}_M(Y) : Y \text{ is a subvariety of } \mathcal{G}_M \}$.

Finally, for every $M'\in \mathcal{M}_n^d$ and $Y\subset \mathcal{G}_M$ a subvariety, there exists a unique $M\in {M}_n^d$ such that $h(\mathbb{V}_{M'}(Y))={h}(\mathbb{V}_{M}(\sigma(Y)))$ for some $\sigma\in \mathfrak{S}_n$. The assignment $\mathbb{V}_{M'}(Y)\mapsto \mathbb{V}_{M}(\sigma(Y))$ is precisely the quotient map $\Gr(d,n)^T\xrightarrow{}\Gr(d,n)^T/\mathfrak{S}_n$, and the result follows.\fin
\end{dem}

By Theorem \ref{V(Y)_2} above, in order to understand the  questions for the map $h$ of~\eqref{Homology_map}, 
it is enough to understand the quotient map $\Gr(d,n)^T\xrightarrow{}\Gr(d,n)^T/\mathfrak{S}_n$ and the map  
$$\theta :\bigsqcup_{M\in {M}_n^d}\{Y : Y \text{ is a subvariety of } \mathcal{G}_M \}\xrightarrow{}H_*(\mb{G}(d,n),\mathbb{Z}),\quad \theta(Y)=h(\mathbb{V}_M(Y))$$
(or on scheme-theoretic points of $\mathcal{G}_M$).

The first one is determined by the study of the fibers of $\mathcal{M}_n^d\xrightarrow{}\mathcal{M}_n^d/\mathfrak{S}_n$, and this amounts to compute the $\mathfrak{S}_n$-orbit $\mathfrak{S}_n\cdot M$ of every $M\in M_n^d$. For this reason we will focus on the map $\theta$, whose image is the set of prime $T$-classes. 

For every prime $T$-class $\lambda \in H_*(\mb{G}(d,n),\mathbb{Z})$ and $M\in {M}_n^d$, we define $\mathcal{G}_M(\lambda)=\{Y\subset \mathcal{G}_M\text{ is a subvariety}\::\:\theta(Y)=\lambda\}$ and  ${M}_n^d(\lambda)=\{M\in {M}_n^d\::\:\mathcal{G}_M(\lambda)\neq\emptyset \}$. Then  we can write  $\theta^{-1}(\lambda)=  \bigsqcup_{M\in {M}_n^d(\lambda)}\mathcal{G}_M(\lambda)$. 

\begin{cor}
\label{cor_param_space}
Let $\lambda \in H_*(\mb{G}(d,n),\mathbb{Z})$ be a prime $T$-class. Then $h^{-1}(\lambda)$ equals  
\begin{equation}
\label{param}
    \mathcal{G}(d,n,\lambda):=\bigsqcup_{M\in {M}_n^d(\lambda)}\mathcal{G}_M(\lambda)^{\sqcup\#(\mathfrak{S}_n\cdot M)}.
\end{equation}
\end{cor}
\begin{dem}
Let $\lambda\in\text{Im}(h)$. Then 
%
the result follows from the definition of $\mathcal{G}_M(\lambda)$ and
the fact that $\mathcal{G}_M(\lambda)\cong \mathcal{G}_{M'}(\lambda)$ if  $M'\in \mathfrak{S}_n\cdot M$.\fin
\end{dem}

This reduces our problem to the computation of the spaces \eqref{param}. 
These spaces will be the subject of many of our subsequent results.
For instance, note that if $y\in \mathcal{G}_M$ is a  point, then $\mathbb{V}_M(y)\subset G_M$ is a single torus orbit closure, which we may write as $\overline{Tx}$ for any closed point $x\in G_M$. 
If we denote by $\lambda(M)$ the homology class of this orbit, we will show in Corollary \ref{cor_char_orb} in the next section that $h^{-1}(\lambda(M))$ is parametrized by $\mathcal{G}_M^{\sqcup(\#\mathfrak{S}_n\cdot M)}$.

\subsection{Classes of orbits}
\label{Subsection_Classes_Of_Orbits}
In this section we show that a $T$-invariant subvariety of~$\Gr(d,n)$ with the homology class of an orbit closure is in fact an orbit closure.
The main idea of the argument is an examination of connected components of matroids.

Let $E_\bullet=(E_1,\ldots,E_k)$ be an ordered set partition of~$[n]$,
and let $d_\bullet=(d_1,\ldots,d_k)$ be a $k$-tuple of nonnegative integers whose sum is~$d$.
Let $\mb{C}^{E_i}$ denote the coordinate subspace of $\mb{C}^n$
spanned by the standard basis vectors $\{e_j:j\in E_i\}$.
The direct sum decomposition $\mb{C}^n=\bigoplus_{i=1}^k\mb{C}^{E_i}$
provides a closed embedding
\begin{equation}
\label{product_embedding}
\iota(E_\bullet,d_\bullet):\prod_{i=1}^k \Gr(d_i,\mb{C}^{E_i})\hookrightarrow \Gr(d,n).
\end{equation}
Note that $\Gr(d,n)$ is a one-point space for $d=0$. 

If $M$ is a matroid of rank~$d$ on $n$ elements,
then $G_M$ lies in the image of $\iota(\comp(M),d_\bullet)$, 
where, if $\comp(M) = \{E_1,\ldots,E_k\}$,
then $d_i=\#(B\cap E_i)$ for any basis $B$ of~$M$.
This is the finest partition $E_\bullet$
such that $G_M$ is in the image of $\iota(E_\bullet,d_\bullet)$.

In \eqref{product_embedding} we identify the integral homology of the source with the tensor product
\[\bigotimes_{i=1}^k H_*(\Gr(d_i,\#E_i),\mb{Z}),\]
exploiting the vanishing of odd-degree homology in the K\"unneth theorem. Then, we have the following result.
We denote by $\delta$ the  isomorphism $H^{2(\dim \Gr(d,n)-p)}(\Gr(d,n),\mathbb{Z})\xrightarrow{\sim} H_{2p}(\Gr(d,n),\mathbb{Z})$ induced by Poincar\'e duality.
\begin{p}\label{p_iota}
The pushforward in homology
\[\iota_*:\bigotimes_{i=1}^k H_*(\Gr(d_i,\#E_i),\mb{Z})\to H_*(\Gr(d,n),\mb{Z})\]
along the inclusion of \eqref{product_embedding} is given by
\[\iota_*(\delta(\lambda_1)\otimes\cdots\otimes\delta(\lambda_k)) = \delta(\lambda_1\cdots\lambda_k).\]
\end{p}
In other words, the inclusion in homology from products of sub-Grassmannians that we get from disconnected matroids is Poincar\'e dual to the product in cohomology. 
\begin{dem}
It is enough to do $k=2$, and show that $\iota_*$ acts correctly on Schubert classes.
For an integer~$e$ and a sequence $a\in\mathbb Z^e$, let 
\[\delta_e(a)=(n-d-a_e, \ldots, n-d-a_1).\]
The dual Schubert variety to $\Sigma_a$ under the intersection pairing on~$\Gr(d,n)$ is~$\Sigma_{\delta_d(a)}$.
Also let $\delta_e(\sigma_a)=\sigma_{\delta_e(a)}$, and extend linearly to linear combinations of Schubert classes.

For $i=1,2$, the $d$-subspaces $L\subseteq\mb{C}^n$ which lie in the image of $\iota=\iota(E_\bullet,d_\bullet)$
are those for which $\dim(L\cap E_i)=d_i$.
Let $\mathcal{V}^i$ be a complete flag in which~$E_i$ appears.
Let $\Sigma^{\mathcal{V}^i}_{a^i}$ be a Schubert variety 
that parametrizes spaces $L$ with $\dim(L\cap E_i)=d_i$
and imposes no conditions on intersections on $L\cap\mathcal{V}^i_j$ for any $j>\dim E_i$,
but may impose any Schubert conditions on the lower intersections.
Then the last $d-d_i=d_{3-i}$ entries of the sequence $a^i$ equal $0$.
The intersection $\Sigma^{\mathcal{V}^1}_{a^1}\cap \Sigma^{\mathcal{V}^2}_{a^2}$ is precisely 
$\iota(\Sigma^{t\mathcal{V}^1}_{ta^1},\Sigma^{t\mathcal{V}^2}_{ta^2})$,
where $t\mathcal{V}^i$ and $ta^i$ are the initial segments of $\mathcal{V}^i$ and $a^i$ of length $\dim E_i$.
For some choice of the flags $\mathcal{V}^i$ (for example, as opposite flags \cite{BrionLakshmibai})
the two Schubert varieties $\Sigma^{\mathcal{V}^i}_{a^i}$ meet transversely.
So their intersection $\iota(\Sigma^{t\mathcal{V}^1}_{ta^1},\Sigma^{t\mathcal{V}^2}_{ta^2})$
has cohomology class $\sigma_{a^1}\,\sigma_{a^2}$
and homology class $\delta_d(\sigma_{a^1}\,\sigma_{a^2})$.  
The proof is completed by the fact that
\[\delta_d(\sigma_{a^1}\,\sigma_{a^2})
= \delta_{d_1}(\sigma_{a^1})\,\delta_{d_2}(\sigma_{a^2})
= \sigma_{\delta_{d_1}(ta^1)}\,\sigma_{\delta_{d_2}(ta^2)}. \hfill\blacksquare\]
\end{dem}

Recall that $\mathcal{M}^d_n(p)$ is the set of matroids $M\in\mathcal{M}^d_n$ such that $\dim\Delta(M)=p$, that is, such that $\#\comp(M)=n-p$.
Let $c_p:\mathcal{M}^d_n\xrightarrow{}H^{2(\dim \Gr(n,d)-p)}(\Gr(d,n),\mathbb{Z})$ be the map sending $M$ to the cohomology class $\delta\circ h(\mb{V}_M(*))=[\mb{V}_M(*)]$ of $\mb{V}_M(*)$, where $*$ represents a closed point of $\mathcal{G}_M$.
\begin{p}
\label{TheoremSF}
The map $c_p$ is an additive matroid valuation.
\end{p}
\begin{dem}
For each natural $p$, let $F_p$ be the subgroup of $K^0(\Gr(d,n))$ generated by the classes of structure sheaves of subvarieties of dimension at most $p$. 
Then $F_p/F_{p-1}$ is isomorphic, up to tensoring with $\mathbb Q$,
to the graded component of degree $\dim \Gr(n,d)-p$ of the Chow ring of $\Gr(d,n)$,
by an isomorphism sending the structure sheaf of any $p$-dimensional subvariety to its cycle class,
and of any lower-dimensional subvariety to~0 \cite[Example 15.2.16]{Fulton}.
The Chow ring of the complex Grassmannian is in turn isomorphic to its singular cohomology ring over $\mathbb Z$.

In \cite{FS}, the authors show that the map $y:\mathcal{M}^d_n\longrightarrow K^0(\Gr(d,n))$ sending $M$ to the $K$-class of $\mathcal{O}_{\mb{V}_M(*)}$ is a matroid valuation. 
For each $p$, the image of $\mathcal{M}^d_n(p)$ under~$y$ lands in~$F_p$.
Let $\hat c_p$ be the map obtained by 
composing the restriction of $y$ to $\bigcup_{q\le p}\mathcal{M}^d_n(q)\to F_p$
with the quotient map to $F_p/F_{p-1}$ and the isomorphisms above.
Since the maps we have composed with are group homomorphisms, $\hat c_p$ is a matroid valuation,
and it takes value zero on $\bigcup_{q<p}\mathcal{M}^d_n(q)$,
so $c_p$ is an additive matroid valuation.
\fin
\end{dem}


For the rest of this section, let us use the abbreviated notation $s_a:=\delta(\sigma_a)\in H_*(\Gr(d,n),\mb{Z})$.
Let $h(d,n-1)$ be the sequence $(n-d,1,\ldots,1)$ with $d$ entries and therefore sum $n-1$.  We will call $s_{h(d,n-1)}$ a \textbf{hook}.
The next proposition is essentially due to Speyer.
See \cite[\S7.3]{White} for the definition and basic properties of the beta invariant of a matroid.

\begin{p}\label{p_beta}
Let $M$ be a connected matroid of rank~$d$ on $n>1$ elements.
Then the coefficient of $s_{h(d,n-1)}$ in $c_{n-1}(M)$ is the beta invariant $\beta(M)$. 
In particular, this coefficient is strictly positive.
\end{p}

\begin{dem}
The case where $M$ is representable over~$\mb{C}$ is \textup{\cite[Thm 5.1]{Speyer}}.
The general case follows from this because both sides of the equality are matroid valuations,
and matroids representable over~$\mb{C}$ span all matroids up to valuation.\fin
\end{dem}

We will detect classes of subvarieties contained in the image of some nontrivial inclusion $\iota(E_\bullet,d_\bullet)$
by inspecting the hook content in their homology classes.
Using Proposition~\ref{p_iota}, this ultimately becomes a product on the algebra of symmetric functions $\Lambda$.
The following lemma is routinely verified by any of the forms of the Littlewood-Richardson rule.

\begin{lem}\label{lem_product_of_hooks}
The $\mb{Z}$-submodule of $\Lambda$ generated by Schur functions $s_a$ with $a_2\geq2$ is an ideal $I$.  We have
\[s_{h(d,n-1)}s_{h(d',n'-1)} \equiv s_{h(d+d'-1,n+n'-2)} + s_{h(d+d',n+n'-2)} \pmod I.\]
\end{lem}

Note that if $a_2\not\geq 2$ then $s_a$ is either a hook or equal to~$1$.

\begin{teo}
\label{teo_orbits}
Let $Y\subset \Gr(d,n)$ be a $T$-invariant subvariety such that
$[Y] = [\overline{T x}]$ for some closed point $x\in \Gr(d,n)$.
Then $Y$ is itself a $T$-orbit closure.
\end{teo}

\begin{dem}
Let $M$ be the matroid such that $\overline{T x}=\mb{V}_M(*)$.
By Theorem~\ref{V(Y)}, a subvariety of $\Gr(d,n)$ of complex dimension~$p$ that is not an orbit
is instead $\mb{V}_N(Y)$ 
for some matroid $N$ and positive-dimensional subvariety $Y\in\mathcal{G}_N$.
The torus fibers in this presentation have dimension $\dim T_N=p-\dim Y<p$.
By Theorem~\ref{thm_geometry}(\ref{dim(T)}),
we have $\#\comp(M)=n-p<n-p+\dim Y=\#\comp(N)$.

The class $[\overline{T x}]$ is what was called $c_p(M)$ above.
By Proposition~\ref{p_iota}, 
$c_p(M)$ is the product (via the Poincar\'e duality with cohomology, or in the symmetric algebra)
of the classes $c_{d_i}(E_i)$, where $\comp(M) = \{E_1,\ldots,E_{n-p}\}$
and $d_i=\#(B\cap E_i)$ for a basis $B$ of~$M$.
If $\#E_i=1$ then the unique element of $E_i$ is either a loop (if $d_i=0$) or a coloop (if $d_i=1$), in which cases $c_{d_i}(E_i) = 1$.
Otherwise, by Proposition~\ref{p_beta}, $c_{d_i}(E_i) = \beta_i s_{h(d_i,\#E_i-1)}$ mod~$I$, where $I$ is the ideal of Lemma~\ref{lem_product_of_hooks} and $\beta_i$ is a positive integer,
because $h(d_i,\#E_i-1)$ is the only hook that indexes a class in $H_{2(\#E_i-1)}(\Gr(d_i,\#E_i),\mb{Z})$.
Without loss of generality we may number the connected components so that $E_1,\ldots,E_q$ are exactly the components that are not loops or coloops.
Using earlier notation, $\#E_1+\cdots+\#E_q=n-\#L(M)-\#C(M)$,
so $\sum_{i=1}^q(\#E_i-1)=n-\#\comp(M)$.
By Lemma~\ref{lem_product_of_hooks}, we have
\[c_p(M) \equiv \Big(\prod_{i=1}^{q}\beta_i\Big)\sum_{i=0}^{q-1}\binom qi s_{h(d_1+\cdots+d_q-i,n-\#\comp(M)))}
\pmod I.\]

We now perform a similar analysis for~$N$.
Let $\comp(N) = \{E'_1,\ldots,E'_k\}$
and $d'_i=|B\cap E'_i|$ for a basis $B$ of~$N$.
If a hook $s_{h(i,j)}$ indexes a class in $H_*(\Gr(d'_i,\#E'_i),\mb{Z})$,
then $j\leq\#E'_i-1$.
Proposition~\ref{p_iota} and Lemma~\ref{lem_product_of_hooks} then show that,
if any hooks appear in $[\mb{V}_N(Y)]$,
they are $s_{h(i,j)}$ where 
$j\leq \sum(\#E'_i-1)=n-\#\comp(N)$.
But this cannot have the right degree: $[\mb{V}_N(Y)]=[c_p(M)]$ should belong to $H_{2(n-\#\comp(M))}$, and $j=n-\#\comp(M)$ is not attainable in the above analysis because $\#\comp(M)<\#\comp(N)$.\fin
%
%
\end{dem}

\begin{obs}
A matroid $M$ is \emph{not} uniquely determined by the class $\theta(y)$ for $y\in\mathcal{G}_M(\mathbb{C})$.
For example, there exist two nonisomorphic sparse paving matroids in $M_6^3$ with two circuit-hyperplanes, both representable over $\mathbb{C}$,
and torus orbits in $\mb{G}(3,6)$ corresponding to these two have the same homology class.
However, when $d=2$, $M$ is uniquely determined by $\theta(y)$: see Corollary~\ref{cor_char_orb}.
\end{obs}

\begin{obs}
According to Theorem \ref{teo_orbits}, we get a dichotomy for prime $T$-classes: either all the (integral, $T$-invariant) subvarieties homologous to $\lambda$ are (closures of) $T$-orbits, or no (closures of a) $T$-orbit can be homologous to $\lambda$. Thus it becomes  important to be capable of computing cohomology classes of torus orbits.
\end{obs}

\begin{cor}
\label{cor_computation_orbits}
 Any torus orbit class in $H^*(\mb{G}(d,n),\mathbb{Z})$ can be expressed as a (finite) sum of torus orbit classes corresponding to rigid matroids.
\end{cor}
\begin{dem}
Given $M\in M_n^d$, if $x\in\mathcal{G}_M$ then $\overline{Tx}=\mathbb{V}_M(*)=\text{Tor}(\Delta_M)$. If $M$ is rigid there is nothing left to do, so we suppose that this is not the case. 
We apply induction on the number of integer points in $\Delta_M$, which is the number of bases of~$M$.
Let $\{\Delta_{M_1},\ldots,\Delta_{M_k}\}$
be any nontrivial matroidal subdivision  of $\Delta_M$.
We apply Proposition \ref{TheoremSF} to find 

$$[\mathbb{V}_M(*)]=\sum_i[\mathbb{V}_{M_i}(*)].$$

Every $\Delta_{M_i}$ has fewer integer points than $\Delta_M$, 
so by the inductive hypothesis we can write the right hand side as a sum of rigid matroid classes.
\fin
\end{dem} 

Note that the expression in $H^*(\mb{G}(d,n),\mathbb{Z})$ is independent  of the choice of the maximal matroidal subdivision, since the map $M\mapsto[\mathbb{V}_M(*)]$ is an additive matroid valuation. 

We will characterize the rigid elements of $M\in M_n^2$ in Corollary \ref{cor_rigid}.  We will show in Section \ref{Section_Computing} how Corollary \ref{cor_computation_orbits} can be effectively used to compute the homology class of any torus orbit when $d=2$.
 
\section{The case of Grassmannians of lines}
\label{Section_case_lines}
In this section we develop some aspects of our general theory for the case $d=2$, as an example to showcase its applications,
culminating in Section~\ref{Sect_geo} on the geometry of thin Schubert cells.
Several of the facts in this section are standard in matroid theory, 
and several others echo Kapranov's work \cite{K} (especially Section~1.3 thereof), but we give proofs for self-containedness.

\label{Subsection_Partitions}
We begin with setup, first of all for integer partitions,
which are basic objects in combinatorics. The reader may consult \cite{RS} for a wealth of information, of which we only need the very rudiments.
The importance of partitions to us, in brief, is that torus orbits in $\Gr(2,n)$ are configurations of $2\le\ell\le n$ points on $\mathbb P^1$, and we need the combinatorial information of the sizes of subsets of these $\ell$ points that coincide.
The same combinatorics can be captured by matroids of rank~$2$; this is the subject of Section~\ref{Sect_M2P}.

\begin{defi}\label{defi:partition}
An (integer) \textbf{partition} is a finite multiset of positive integers.
Partitions are normally written as sorted tuples $\pi=(k_1,\ldots,k_\ell)$
where $k_1\ge k_2\ge\cdots\ge k_\ell\geq1$.
We call the $k_i$ the \textbf{parts} of $\pi$. 
We define two statistics $w,\ell$ on partitions:
$w(\pi)=k_1+\cdots+k_\ell$ is the \textbf{weight} of the partition $\pi$ and 
$\ell(\pi)=\ell$ its \textbf{length}.
\end{defi}
That is, $\ell(\pi)$ is the number of parts of $\pi$, and $\pi$ is a partition of the integer $w(\pi)$.
As a purely formal shorthand, we may also use superscripts to write partitions with identical parts collapsed. 
If we collapse every set of identical parts then we have $\pi=(k_{i_1}^{\alpha_1},\ldots,k_{i_t}^{\alpha_t})$, where  $k_{i_1}>\cdots>k_{i_t}$ and $\alpha_j\geq1$ is the \textbf{multiplicity} of the part $k_{i_j}$.

\begin{defi}
\label{def_partition}
The set $\Pi_n$ of \textbf{partial partitions} of~$n$ is the set of partitions of length at least two and weight at most $n$.
\end{defi}

The weight and length statistics satisfy $2\leq \ell(\pi) \leq w(\pi)\leq n$ for any $\pi\in\Pi_n$.

%
%
%
%

\subsection{Matroids of rank 2 as partitions}
\label{Sect_M2P}
Matroids of rank~$2$ have an easy structure.
They are entirely determined by their set of loops and their partition into connected components.
Note that this is a set partition, not one of the integer partitions of the previous section,
although later we will ``forget the labelling'' of the ground set and pass to integer partitions.

\begin{p}\label{p:basis_partition}
For a fixed natural number $n$,
matroids $M$ of rank~$2$ on ground set $[n]$ are in bijection with
families $\pi=\pi(M)=\{P_1,\ldots,P_l\}$ of $l\geq2$ pairwise disjoint non-empty sets $P_i\subset [n]$.
The bases of~$M$ are the \textbf{partial 2-transversals} of~$\pi$,
that is the sets $\{i_1,i_2\}$ with $i_k\in P_{a_k}$ for distinct elements $P_{a_1},P_{a_2}$ of~$\pi$.
\end{p}

In the spirit of cryptomorphism, we will often specify $M$ by giving the data $([n],\pi)$.
In the sequel we will use the notation $S(M):=\bigcup_{i=1}^l P_i$ for the set of non-loops of~$M$.

\begin{dem}
Given a matroid $M=([n],\mathcal B(M))$ of rank 2, we will show how to associate to it a partition $\pi$ of a subset of $[n]$. One can check that the association is bijective.
Let $S$ be the set of elements of~$M$ which are not loops, i.e.\ are contained in some basis. 
Let $R\subset S\times S$ be the relation for which $iRj$ if and only if $\{i,j\}$ is not a basis of~$M$, so in particular $iRj$ if $i=j$. 
This is an equivalence relation on~$S$, known as the \textbf{parallelism} relation:
if $\{a,c\}\in\mathcal B(M)$ and $b\in S$, applying \eqref{mep} with $I=\{b,i\}\in\mathcal B(M)$ and $J=\{a,c\}$ shows that $(a,b,c)$ is not a counterexample to transitivity. 
Take $\pi$ to be the  partition that $R$ induces on~$S$. Then $M=([n],\pi)$. 
\fin
\end{dem}

Our parameter functions $(w,\ell)$ for partitions take the following form for a matroid $M\in \mathcal{M}_n^2$: we have $w(M):=n-\#L(M)$, and $\ell(M)$ is the number of connected components of the loopless matroid $M\setminus L(M)$. 

Labelling matroids in $\mathcal{M}_n^2$ by set partitions of~$[n]$ gives us a simple description of the orbit set $\mathcal{M}_n^2/\mathfrak{S}_n$: orbits are labelled by $\Pi_n$, the set of partial partitions of the integer $n$ with at least two parts.
We also construct a concrete set of orbit representatives $M_n^2=\{M_\pi\in \mathcal{M}_n^2\::\:\pi\in\Pi_n\}$ indexed by~$\Pi_n$. 

For $\pi=(k_1,\ldots,k_\ell)\in\Pi_n$ with $k_1\geq\cdots\geq k_\ell$ and $w(\pi)=\sum_ik_i$\,, we will form a particular set partition 
$\setpart\pi$ of $[w(\pi)]=\{1,\ldots,w(\pi)\}$,
whose parts are the sets of consecutive integers
\begin{equation}
\label{particular_partition}
P_i=\{k_{i+1}+\cdots+k_{\ell}+1,\, \ldots,\, k_{i}+\cdots+k_{\ell}\}.
\end{equation}
for $1\le i\le\ell$. 

\begin{defi}\label{def_M_pi}
We denote by $M_\pi$ the matroid $([n],\setpart\pi)$ of rank~2 defined on $[n]$ by the set partition $\setpart\pi$ as in Proposition~\ref{p:basis_partition}.  
\end{defi}

So $M_\pi=M_\pi^{\rm loopless}\oplus (U_1^0)^{\oplus {[n-w(\pi)]}}$, where $M_\pi^{\rm loopless}=([w(\pi)],\setpart\pi)$ is the loopless matroid of rank~2 defined by $\setpart\pi$, and $(U_1^0)^{\oplus {[n-w(\pi)]}}$ is the uniform matroid on $[n-w(\pi)]$ of rank 0, i.e.\ consisting only of loops.

\begin{exa}
Let us consider $\pi=(3^2,2^3,1)\in\Pi_{16}$. Then the corresponding partition of  $[w(\pi)]=[13]$ is  
\[\setpart\pi=\big\{P_6=\{1\},P_5=\{2,3\},P_4=\{4,5\},P_3=\{6,7\},P_2=\{8,9,10\},P_1=\{11,12,13\}\big\},\]
and $\mathcal B(\setpart\pi)$ consists of the $69$ possible pairs $\{\{i,j\}\::\:1\leq i< j\leq 13\}$ where $i$ and $j$ belong to different sets of the above partition. The elements of $[16]$ not in a part of $\setpart\pi$ are the loops of $M_\pi=([16],\mathcal B(\setpart\pi))$, namely $\{14,15,16\}$.
\end{exa}  

The reason we have chosen a ``backward'' set partition with its parts in order of increasing size
is to arrange that 
$M_\pi$ is a Schubert matroid when $\pi$ is of the form $(k,1^m)$.
This will be convenient in Section~\ref{Subsection_Thin_SC} when we label Schubert varieties in $\Gr(2,n)$ with such partitions.

\begin{teo}
\label{theo_complete_set}
A complete set of representatives $M_n^2$ for $\mathcal{M}_n^2/\mathfrak{S}_n$ is 
\begin{equation*}
 \{M_\pi
\::\:\pi\in\Pi_n\}.
\end{equation*} 
\end{teo}
\begin{dem}

Pick an arbitrary $M\in \mathcal{M}_n^2$. Under the action of $\mathfrak{S}_n$, we can suppose that the loops of $M$ are $\{n-\#L(M)+1,n-\#L(M)+2,\ldots,n\}$, so that the restriction $M|[n-\#L(M)]$ is a loopless matroid. By applying the previous construction, we find a partition $P=\{P_1,\ldots,P_l\}$ of  $[n-\#L(M)]$ such that  $M\cong M_\pi$ for $\pi=\#P_1+\cdots+\#P_l$. Finally, it is clear that different partitions represent different $\mathfrak{S}_n$-orbits.\fin
\end{dem}

So the projection $\mathcal{M}_n^2\longrightarrow\mathcal{M}_n^2/\mathfrak{S}_n$ allows us to go from a set partition of some $S\subset [n]$ with $\#S\geq2$ to a partition of the integer $2\leq\#S\leq n$. The assignment $\pi\mapsto M_\pi$ gives a section of this projection.

Now that we have a complete set of representatives $M_n^2$ for $\mathcal{M}_n^2/\mathfrak{S}_n$, we  discuss the number of symmetries in $\mathcal{M}_n^2$ of each representative $M_\pi$. 
For any $\pi=(k_1^{\alpha_i},\ldots,k_t^{\alpha_t})\in \Pi_n$, on this page we denote $\tau_n(\pi)=\#(\mathfrak{S}_n\cdot M_\pi)$. 
\begin{p}
\label{Number_Symmetries}
We have 
\begin{equation*}
\tau_n(\pi)=\frac{n!}{(n-w(\pi))!(k_1!)^{\alpha_1}\alpha_1!\cdots(k_t!)^{\alpha_t}\alpha_t!}.
\end{equation*}
\end{p}
\begin{dem}

We have seen that matroids are pairs $([n],\pi)$, where $\pi$ is a partition of a subset $S\subset [n]$. So, if $S=[n]$, i.e.\ $ w(\pi)=n$, then we just need to count the number of partitions of the set $[n]$  into $l$ subsets of which $\alpha_i$ contain $k_i$ elements, for each unique part $k_1>\cdots>k_t\geq1$. This number is known as an unordered multinomial coefficient:
$$c_\pi:=\frac{n!}{(k_1!)^{\alpha_1}\alpha_1!\cdots(k_t!)^{\alpha_t}\alpha_t!}.$$
This solves the loopless case, because $\tau_{n}(\pi)=c_\pi$. If $S\neq[n]$, that is, $ w(\pi)<n$, then   the numbers $\tau_n(\pi)$ are just multiples of $\tau_{w(\pi)}(\pi)$, namely $\tau_n(\pi)=\binom{n}{w(\pi)}\tau_{w(\pi)}(\pi)$, since any family of $n-w(\pi)$ elements may appear as the set of loops of a matroid.\fin
\end{dem}
\begin{exa}
For Schubert partitions (see p.~\pageref{SchubertPartition}) we have:
\begin{enumerate}
\item  $\tau_n((1^m))=\binom{n}{m}$, and 
\item $\tau_n((k,1^{m-k}))=\binom{n}{m}\binom{m}{k}$, if  $2\leq k<m\leq n$.
\end{enumerate}
\end{exa}
As another example we have that 
when $M_\pi$ is rigid,
the set $\mathcal{G}(2,n,h(\overline{G}_\pi))$ of Corollary~\ref{cor_param_space}
is a finite set of points of cardinality $\#(\mathfrak{S}_n\cdot M_\pi)$ by Corollary~\ref{cor_char_orb}, and this number can be computed from the formula shown in Proposition \ref{Number_Symmetries}.

Through the bijections $\mathcal{M}_n/\mathfrak{S}_n\longleftrightarrow \{M_\pi\::\pi\in \Pi_n\}\longleftrightarrow\Pi_n$, we can associate properties of the matroid $M_\pi$ with the partition $\pi$. 
In this way we define the polytope of $\pi\in\Pi_n$ as $\Delta_\pi=\Delta(M_\pi)$,
and restate Remark~\ref{dimension_of_poly} using partitions.
\begin{cor}
\label{dimension_of_poly_from_partition}
 We have 
 \begin{equation*}
 \text{dim }\Delta_\pi=
 \begin{cases}
  w(\pi)-2,&\text{ if }\ell(\pi)=2,\\
  w(\pi)-1,&\text{ if }\ell(\pi)>2.\\
 \end{cases}
 \end{equation*}
\end{cor}
\begin{dem}

We use the fact that the dimension of $\Delta_{\pi}$ is equal to $\text{dim }\Delta(M_\pi)=n-\#\comp(M_\pi)=w(\pi)-\#\comp(M_\pi)$ for the matroid $M_\pi$.
Suppose first that $\ell(\pi)=2.$ If $w(\pi)=2$, then $\#C(M_\pi)=2$, and the dimension is $2-2=0$. Suppose $w(\pi)>2$. If $\pi=(w(\pi)-1,1)$, then $\#C(M_\pi)=1$, so $M_\pi\cong U_{1}^1\oplus U_{w(\pi)-1}^1$ and the dimension is $w(\pi)-1-1$. If $\pi=(w(\pi)-k,k)$ for $2\leq k\leq w(\pi)-2$, then $M_\pi$ has no co-loops and has two connected components, so it has dimension $w(\pi)-2$. If $\ell(\pi)>2$, then the matroid has no co-loops and it is connected. Then the dimension is $w(\pi)-1$.\fin
\end{dem}

\subsection{The geometry of the thin Schubert cells}
\label{Sect_geo}
We first describe $G_M$ in terms of the matrices whose rows are bases for linear spaces $x\in G_M$
(so-called \textbf{representations} of~$M$).
Let $x_A\in\Gr(2,n)$ be the span of the rows of the rank~$2$ matrix $A\in\text{Mat}_{2\times n}(\mathbb{C})$. We denote the columns of~$A$ as $A[i]$, $i=1,\ldots,n$.
Recall the partition $\pi(M)$ of the set $S(M)$ from Proposition~\ref{p:basis_partition}.
\begin{lem}\label{lem_representations}
We have $x_A\in G_M$ if and only if
\begin{enumerate}
    \item $A[i]=0$ exactly when $i\not\in S(M)$;
    \item for $i,j\in S(M)$, $A[i]$ and $A[j]$ are parallel exactly when $i$ and $j$ are in the same part of $\pi(M)$.
\end{enumerate}
\end{lem}

\begin{dem}
Letting $A[i,j]$ denote the $2\times 2$ matrix with columns $A[i]$ and $A[j]$, the lemma is clear from the fact
\[\mathcal B(M) = \{\{i,j\}\in{\textstyle\binom{[n]}2} : \det A[i,j]\ne0\}.
\hfill\blacksquare\]
\end{dem}


Theorem \ref{theo_complete_set} provides us a complete set of representatives $M_n^2=\{M_\pi\::\:\pi\in\Pi_n\}$ for the space of orbits $\mathcal{M}_n^2/\mathfrak{S}_n$ indexed by a particular set of partitions of the integer $n$.  
We write $G_\pi:=G_{M_\pi}$, and end up with a family 
$\mathcal{B}:=\{T_\pi\xrightarrow{}G_\pi\xrightarrow{}\mathcal{G}_\pi:\:\pi\in M_n^2\}$ 
which we see (modulo the action of $\mathfrak{S}_n$) as a decomposition of $\Gr(2,n)$, where  $T_\pi\xrightarrow{}G_\pi\xrightarrow{}\mathcal{G}_\pi$ denotes the bundle \eqref{intro_fb} corresponding to the matroid $M_\pi\in M_n^2$. 
We explicitly describe these bundles in the next theorem.

We make use of~$\mathcal{M}_{0,n}$, the  moduli space of ordered configurations of $n\ge3$  distinct points on a rational curve up to isomorphism. 
A general reference for $\mathcal{M}_{0,n}$ is~\cite[Section 0]{KV}. 
$\mathcal{M}_{0,n}$ is isomorphic to the configuration space $F(\mathbb{CP}^1\setminus\{0,1,\infty\},n-3)$ of ordered configurations of $n-3$  distinct points in $\mathbb{CP}^1\setminus\{\infty,0,1\}$,
as one can identify the rational curve with $\mathbb{CP}^1$ and then use its 3-transitive automorphism group to bring three of the marked points to the positions $\infty$, $0$, and $1$.
The latter configuration space is in fact the complement of an arrangement of hyperplanes in $\mathbb{P}^{n-3}$, 
namely the hyperplane at infinity,
the $2(n-3)$ hyperplanes $\{x_i=0\}$ and $\{x_i=1\}$,
and the $\binom{n-3}2$ hyperplanes $\{x_i=x_j\}$. 
The next theorem in the case $\pi=(1^n)$ is Kapranov's \cite[Theorem 4.1.8]{K}.

\begin{teo}
\label{thm_geometry}
Consider the bundle $T_{\pi}\xrightarrow{}G_{\pi}\xrightarrow{}\mathcal{G}_{\pi}$ for $M_\pi\in M_n^2$.
\begin{enumerate}
\item The parameter space $\mathcal{G}_{\pi}$ is a point if $\ell(\pi)= 2$, and it is $\mathcal{M}_{0,\ell(\pi)}$ otherwise.
\item\label{dim(T)} The dimension of the torus $T_{\pi}$ is $n-\#\comp(M_\pi)$.
That is, it is $w(\pi)-2$ if $\ell(\pi)=2$, and $w(\pi)-1$ otherwise.
\item The dimension of the thin Schubert cell $G_{\pi}$ is $w(\pi)+\ell(\pi)-4$.
\end{enumerate} 
\end{teo}

\begin{dem}
A $T$-orbit of matrices $A\in\text{Mat}_{2\times n}(\mathbb{C})$, seen columnwise, is an $n$-tuple of $\mathbb{C}^*$ orbits on $\mathbb C^2$, i.e.\ an $n$-tuple $(p_1,\ldots,p_n)$ of elements of $\mathbb P^1\cup\{(0,0)\}$.
For $x_A\in G_\pi$, Lemma~\ref{lem_representations} says that $p_i=(0,0)$ exactly when $i\not\in S(M_\pi)$,
and for $i,j\in S(M_\pi)$, we have $p_i=p_j\in\mathbb P^1$ if and only if $i$ and~$j$ are in the same part of~$\pi$.
That is, $T$-orbits of such $A$ are in bijection with choices of a distinct point of $\mathbb P^1$ for each of the $\ell(\pi)$ parts of~$\pi$.
The left $\mathrm{GL}_2(\mathbb{C})$-action on~$A$ induces the action on $\mathbb P^1$ by its automorphisms.
So $\mathcal{G}_\pi=G_\pi/T_\pi=G_\pi/T$ is $\mathcal{M}_{0,\ell(\pi)}$ by definition, 
except when $\ell(\pi)=2$ (and $\mathcal{M}_{0,\ell(\pi)}$ is not defined) in which case $\mathcal{G}_\pi$ is a point because the automorphisms of $\mathbb P^1$ are at least 2-transitive.
This is the first statement.

For the second statement, the dimension of $T_{\pi}$ equals the dimension of $\Delta_\pi$, and this is given by Corollary~\ref{dimension_of_poly_from_partition}. The last part comes from the fact that the dimension of~$G_{\pi}$ is the sum of
$\dim\mathcal{G}_{\pi}=\max\{\ell(\pi)-3,0\}$ and $\dim T_{\pi}$.\fin
\end{dem}

For later use we record the dimension of a general subvariety $\mathbb{V}_{\pi}(y)$.

\begin{cor}\label{p:dim V_y}
The dimension of~$\mb{V}_\pi(y)$ is
\begin{equation*}
\text{dim}_\mb{C}\mb{V}_\pi(y)=
\begin{cases}
w(\pi)+\ell(\pi)-4,&\text{if }\ell(\pi)\leq3,\\
w(\pi)+\ell(\pi)-4- \codimC y,&\text{if }\ell(\pi)\geq4.
\end{cases}
\end{equation*}
\end{cor}

We get a characterization of the rigid elements of $M_n^2$. 
If $M_\pi$ is rigid for some $\pi\in\Pi_n$, we will allow ourselves to say that $\pi$ is rigid.

\begin{cor}
\label{cor_rigid}
A partition $\pi\in \Pi_n$ is rigid if and only if $\ell(\pi)\leq3$.
\end{cor}

\begin{dem}
The ``only if'' direction is the folklore theorem of Lafforgue mentioned after Definition~\ref{def_rigid_matroid}.
If $\ell(\pi)\leq3$ then all connected components of $M_\pi$ are series-parallel so the ``if'' direction is Proposition~\ref{p:sp rigid}.
\end{dem}

We next follow up on the special situation when $d=2$ in Remark~\ref{obs_tscdesc}.
\begin{p}
\label{P_stratification}
The thin Schubert decomposition of $\Gr(2,n)$ (Proposition~\ref{P_stratification_1}) is a topological stratification, i.e.\ the closure of each stratum is a union of strata.
\end{p}
Since the strata are labelled by the matroids in the set $\mathcal{M}_n^2$,
the proposition induces a partial order on~$\mathcal{M}_n^2$.
This is the restriction to $\mathcal{M}_n^2$ of the \textbf{weak order} from matroid theory,
which is defined on the set of all matroids on ground set $[n]$.
\begin{dem}
Let $M\in \mathcal{M}_n^2$. 
Let $\overline{G_M}\subset\Gr(2,n)$ be the closure of $G_M$. We show that $\overline{G_M}$ is a union of strata.
The set $\{x_I\::\:I\in \mathcal B(M)\}$ is a set of projective coordinates for $\overline{G_M}$,
and $G_M$ is the locus where none of these coordinates vanish.

Take a point $x\in\overline{G_M}$, and let its matroid be $M'$.
Then $\mathcal B(M')\subseteq\mathcal B(M)$.
If $\{i,j\}\not\in\mathcal B(M')$, then by definition of the parallelism relation (see the proof of Proposition~\ref{p:basis_partition})
either $i$ and $j$ are parallel nonloops in $M'$ or at least one of them is a loop.
Using this observation iteratively on the elements of $\mathcal B(M)\setminus\mathcal B(M')$,
we see that $M'$ can be obtained from~$M$ by a sequence of steps of these two kinds:
\begin{enumerate}
\item\label{item_stratification_F} we turn an element $i$ into a loop, removing it from $S(M)$;
\item\label{item_stratification_B} if $S(M)$ has at least 3 parts, we diminish the number of parts by one by merging two of them, $P$ and $Q$, together. 
\end{enumerate}

We use the notation $F_i^-(M)$ for the matroid obtained by step~\eqref{item_stratification_F}
and $B^-_{P,Q}(M)$ for the matroid obtained by step~\eqref{item_stratification_B}.
Matroids $N$ obtained from~$M$ iterated application of the operations $F^-$ and $B^-$ are therefore the only candidates for matroids with $G_N\subset\overline{G_M}$.

By Lemma~\ref{lem_representations}, all of these candidate strata do appear in $\overline{G_M}$, and the theorem quickly follows.
The closure of any one-dimensional $\mathbb{C}^*$-orbit $O\subset\mathbb{C}^2$ is $O\cup\{(0,0)\}$;
since the $i$\/th column can vary independently of the remaining columns of~$A$, 
the whole stratum $G_{F^-_i(M))}$ is a subset of $\overline{G_M}$.
Similarly, $G_{B^-_{P,Q}(M)}$ is a subset of $\overline{G_M}$
because in $(\mathbb P^1)^{\ell(pi)}$, seen as a configuration space of points on~$\mathbb P^1$ that are not necessarily distinct,
the set of configurations where only points $P$ and $Q$ coincide lies in the closure of the set of configurations where all points are distinct.
\fin
\end{dem}


\begin{obs} With reference to Corollary~\ref{p:dim V_y},
it follows from Theorem \ref{thm_geometry} that applying the operation $F^-_i$ 
to $M_\pi$ leads to a decrement by 1 in the dimension of the fiber $T_\pi$ of the bundle. 
On the other hand, applying the operation $B^-_{P,Q}$ 
to $M_\pi$ leads to a decrement by 1 in the dimension of the base $\mathcal{G}_\pi$ of the bundle.
We chose the letters $F$ and $B$ to stand for ``fiber'' and ``base''.
\end{obs}

\section{The ring $H^*(\Gr(2,n),\mathbb{Z})$ through partitions}
\label{Section_Ring_part}
The purpose of this section is to express the product structure of $H^*(\mb{G}(2,n),\mb{Z})$ in terms of partitions, for effective computation of our algorithms. We also compute in Theorem \ref{thm_map} the classes of $\overline{G}_\pi$ for every $\pi\in \Pi_n$.
We start with the case of Schubert classes, as a warm up.

\subsection{Schubert partitions}
\label{Section_Schubert}
First we introduce the concept of Schubert partition. 
These will be the partitions indexing Schubert matroids as defined in Section~\ref{Subsection_Decomposition};
for the details see Theorem \ref{EHM}.

\begin{defi}
\label{SchubertPartition}
A partition $\pi\in \Pi_n$ is \textbf{Schubert} if it is of the form $\pi=(k,1^{m-k})$ for $1\leq k< m\leq n$.
\end{defi}

Observe that the case $(1^m)=(1,1^{m-1})$ where $k=1$ is included.

Recall that in Section~\ref{Subsection_Decomposition} we used $\mathcal{P}(a_1,a_2) =(n-1-a_1, n-a_{2})$ to index Schubert varieties. 
A different indexing suits our partitions better. 
Schubert partitions are in bijection with the set $\{(k,m)\::\:1\leq k<m\leq n\}$ where the pair $(k,m)$ is sent to the partition $(k,1^{m-k})$. 
The proof of the next proposition is omitted as elementary. 

\begin{p}
\label{Name_Of_SP}
Let $A_n=\{(a_1,a_2)\in\mathbb{Z}^2\::\:n-2\geq a_1\geq a_2\geq0\}$ as in \eqref{indexing_sequence}. The map $\mathcal{R}:A_n\longrightarrow \Pi_n$ defined as $\mathcal{R}(a_1,a_2)=((a_1-a_2+1),1^{n-a_1-1})$ is an injection, and $\mathcal{R}(A_n)$ is the set of Schubert partitions.
\end{p}
For a Schubert partition $\pi$, we have $\mathcal{R}^{-1}(\pi)=(a_1=n-\ell(\pi),a_2=n-w(\pi))$, and therefore the sum of the coordinates of $\mathcal{R}^{-1}(\pi)$ is $2n-\ell(\pi)-w(\pi)=2n-4-[w(\pi)+\ell(\pi)-4]=\codimC G_\pi$.

\subsection{Thin Schubert classes}
\label{Subsection_Thin_SC}
Thin Schubert cells in $\Gr(d,n)$ can be constructed as intersections of Schubert cells for different flags
arising from different orderings of the standard basis of $\mathbb{C}^n$.
These intersections are highly non-transverse in general, so they do not provide a way to write the cohomology class of the closure of a thin Schubert cell as a product of Schubert classes.
The main result of this section, Theorem~\ref{thm_map}, shows that the situation is nicer when $d=2$:
we will write the cohomology class of $\overline{G}_\pi$ for $\pi=(k_1,\ldots,k_\ell)$ as a product of Schubert classes, one for each part $k_i$ of~$\pi$.
The Schubert classes in question are those in Definition~\ref{defi:sigma_k(m)}.

(Non-thin) Schubert varieties, whose cohomology classes were introduced in Section~\ref{Subsection_Decomposition}, are in fact closures of certain thin Schubert cells.  We handle these first.
Let $\mathcal{V}$ be the standard complete flag in $\mathbb{C}^n$,
that is, $\mathcal{V}=(\{0\}=V_0\subset V_1\subset\cdots\subset V_n=\mathbb{C}^n)$, where $V_i=\mathbb{C}^i=\text{Spec}(\mathbb{C}[x_j\::\:1\leq j\leq i])$.
Recalling the definitions of Section~\ref{Subsection_Decomposition},
$M\in \mathcal{M}_n^2$ is a Schubert matroid $M(\Sigma^\mathcal{V}_{a_1,a_2})$ if there exists a Schubert variety $\Sigma^\mathcal{V}_{a_1,a_2}\subset \Gr (2,n)$ such that $M$ is the matroid associated to generic elements of $\Sigma^\mathcal{V}_{a_1,a_2}$. 
Recall the map $\mathcal{R}$ from Proposition \ref{Name_Of_SP}.

\begin{defi}\label{defi:sigma_k(m)}
Define $\Sigma_k(m):=\Sigma^\mathcal{V}_{\mathcal{R}(a_1,a_2)}$ and 
$\sigma_k(m):=[\Sigma_k(m)]=\sigma_{\mathcal{R}(a_1,a_2)}$,
where $k=a_1-a_2+1$ and $m-k=n-a_1-1$.
\end{defi}
That is, $\mathcal{R}(a_1,a_2)=(k,1^{m-k})$. 
Naming this partition $\pi=(k,1^{m-k})$, the matroid $M(\Sigma_k(m))$ equals $M_\pi$.
Our canonical representative of the integer partition $\pi$, as in \eqref{particular_partition},
is the set partition $\pi=\{I_1,\ldots,I_{m-k+1}\}$ of $[m]$ defined as 
$I_j=\{j\}$ for $j=1,\ldots,m-k$ and $I_{m-k+1}=\{m-k+1,\ldots,m\}$.
We see that $M(\Sigma_k(m))$ 
is the matroid of partial 2-transversals of~$\pi$.
Notated using the parameters $k$ and $m$, equation \eqref{eq_general_Schubert_matroid} becomes the following.

\begin{teo}[Theorem 4.3 from \cite{EH}]
\label{EHM}
Let $(k,1^{m-k})=\mathcal{R}(a_1,a_2)$ and  set $\Sigma_k(m):=\Sigma^\mathcal{V}_{\mathcal{R}(a_1,a_2)}$. The set of bases $B_k(m)$ of $M(\Sigma_k(m))$ is
\[
B_k(m)=\left\{I=(i_1<i_{2})\in\binom{[n]}{2}\::\:i_1\leq m-k, i_2\leq m\right\}.
\]
\end{teo}
%
The Schubert variety $\Sigma_k(m)$ equals the closed thin Schubert cell $\overline{G}_M$ where $M=M(\Sigma_k(m))$.
Therefore the cohomology class $[\overline{G}_M]$ equals $\sigma_k(m)$.


\begin{obs}
We have
$$\Sigma_{k}(m)=\{W\in \Gr (2,n)\::\:W\subset V_m,\:W\cap V_{m-k}\neq\{0\}\}=\{W\in \mb{G}(2,m)\::\:W\cap V_{m-k}\neq\{0\}\}.$$
Therefore our coordinate $m$ expresses the position of~$W$ in the filtration of Grassmannians $\{*\}=\mb{G}(2,2)\subset \mb{G}(2,3)\subset \cdots\subset \Gr (2,n)$
induced by the inclusions among the $\mathbb C^i=\operatorname{span}\{e_1,\ldots,e_i\}$.
\end{obs}

In preparation for the next proof, we describe the expansion of the product of two Schubert classes for $d=2$ back into Schubert classes, both in standard notation and in the notation of partitions.
Write $|(a_1,a_2)|=a_1+a_2=\codimC\Sigma^\mathcal{V}_{a_1,a_2}$. 
Assuming $a_1-a_2\geq b_1-b_2$ without loss of generality, the formula for the product is 
\begin{equation}\label{computing_product}
\sigma_{a_1,a_2}\,\sigma_{b_1,b_2}=\sum_{\substack{|c|=|a|+|b|\\a_1+b_1\geq c_1\geq a_1+b_2}}\sigma_{c_1,c_2},
\end{equation}
where if $c_1>n-2$ the symbol $\sigma_{c_1,c_2}$, which is not a Schubert class in $\Gr (2,n)$, is taken to equal~0. See \cite[Proposition 4.11]{EH}. 
In terms of Schubert partitions $\pi_1,\pi_2$, \eqref{computing_product} reads:

$$\gamma(\pi_1)\cdot\gamma(\pi_2)=\sum_{\substack{|\pi_3|=|\pi_1|+|\pi_2|\\2n-\ell(\pi_1)-\ell(\pi_2)\geq n-\ell(\pi_3)\geq 2n-\ell(\pi_1)-w(\pi_2)}}\gamma(\pi_3),$$
where $|\pi|=2n-\ell(\pi)-w(\pi)$ and we are assuming that $w(\pi_1)-\ell(\pi_1)\geq w(\pi_2)-\ell(\pi_2)$.





We now construct an injective map $\gamma:\Pi_n\longrightarrow H^*(\Gr (2,n),\mathbb{Z})$ such that $\gamma(\pi)$ is the cohomology class of $\overline{G}_\pi$,
agreeing with the association above when $\overline{G}_\pi$ is a Schubert variety.
We call then $\gamma(\Pi_n)$ the set of thin Schubert classes for obvious reasons.

\begin{teo}
\label{thm_map}
Consider the map $\gamma:\Pi_n\longrightarrow  H^{*}(\Gr (2,n),\mathbb{Z})$ defined on $\pi=(k_1,\ldots,k_\ell)$ as
\begin{equation}\label{eq:def gamma}
\gamma(\pi):=\sigma_1(n-1)^{n-w(\pi)}\cdot\prod_{i=1}^\ell \sigma_{k_i}(n).
\end{equation}
Then 
the cohomology class of $\overline{G}_\pi$   is $\gamma(\pi)$.
Moreover, $\gamma$ is injective.
\end{teo}
\begin{dem}

To prove that $\gamma(\pi)=[\overline{G}_\pi]$, we express $\overline{G}_\pi$ as a transverse intersection of Schubert varieties in a way matching the right side of~\eqref{eq:def gamma}.

An atlas of affine charts for $\Gr(2,n)$ is given by the nonvanishing loci of each individual homogeneous coordinate $\{x_I\::\:I\in \binom{[n]}{2}\}$ of the Grassmann-Pl\"ucker embedding.
On the chart indexed by $I=\{i_1,i_2\}$, the matrices $A$ of Lemma~\ref{lem_representations} can be chosen so that $A[i_1,i_2]$ is the $2\times 2$ identity matrix; with this choice, the other entries of~$A$ are a system of affine coordinates.
Choose a chart which meets $\overline{G}_\pi$, which is to say $I\in\mathcal B(M_\pi)$.
By Lemma~\ref{lem_representations}, 
the equations of $\overline{G}_\pi$ in this chart are $A[i]=0$ for $i\not\in S(M_\pi)$,
and $\det A[i,j]=0$ for $i,j$ in the same part of $\pi(M_\pi)=\setpart\pi$.

Now observe that for each individual $i\not\in S(M_\pi)$, the two equations $A[i]=0$ cut out $\overline{G}_M$ for the matroid $M=([n],\{\{j\}:j\in[n],j\ne i\})$, which has $i$ as its only loop and no pairs of parallel elements.
Likewise the equations $\det A[j,k]=0$ for $j,k$ both in a single part $P_i$ of~$\setpart\pi$ cut out $\overline{G}_M$ for the matroid $M=([n],P_i\cup\{\{j\}:j\in[n]\setminus P_i\})$, which has no loops and $P_i$ as only nontrivial parallel class.
Since these collections of equations are in disjoint sets of variables in every chart, the various $\overline{G}_M$ just listed intersect transversely, and their intersection is~$\overline{G}_\pi$.
Therefore $[\overline{G}_\pi]$ is the product of these various $[\overline{G}_M]$.
The proof finishes by identifying $\overline{G}_M$ for $M=([n],\{\{j\}:j\in[n],j\ne i\})$ as the Schubert variety $\Sigma_1(n-1)$
and $\overline{G}_M$ for $M=([n],P_i\cup\{\{j\}:j\in[n]\setminus P_i\})$ as the Schubert variety $\Sigma_{k_i}(n)$, where $k_i=|P_i|$,
with respect to suitable flags in the orbit $\mathfrak{S}_n\cdot\mathcal{V}$.
%

To prove that $\gamma$ is injective,
we use the fact that Schubert classes in $H^{*}(\Gr (2,n),\mathbb{Z})$ have the same multiplicative structure constants as Schur functions in two variables,
as discussed around Lemma~\ref{lem_product_of_hooks}.
We spell this out concretely.
If $a_1\ge a_2\ge 0$ are integers, define the polynomial
\[s_{a_1,a_2}=\sum_{i=a_2}^{a_1} x^{a_1+a_2-i} y^i\]
in indeterminates $x,y$. 
The reader may check that in $\mathbb Z[x,y]$ these polynomials are linearly independent, and
the product $s_{a_1,a_2}\,s_{b_1,b_2}$ is the sum of terms $s_{c_1,c_2}$ 
where the indices $(c_1,c_2)$ appearing are precisely those in the sum \eqref{computing_product},
except that we don't set any of the terms to zero.
Temporarily let $\Lambda$ be the ring generated as a free abelian group by all symbols $\sigma_{a_1,a_2}$ with $a_1\ge a_2\ge 0$, with the product specified by~\eqref{computing_product} without setting terms to zero.
Then $\widehat s(\sigma_{a_1,a_2})=s_{a_1,a_2}$ defines an injective ring map
$\widehat s:\Lambda\to\mathbb Z[x,y]$.
It is $\widehat s$ that we use in this proof, but for completeness we note that $\widehat s$ descends to
\[s:H^{*}(\Gr (2,n),\mathbb{Z})\to\mathbb Z[x,y]/\langle s_{c_1,c_2} : c_1>n-2\rangle.\]

Now given $\pi\in\Pi_n$, rewrite the right hand side of \eqref{eq:def gamma}
in the Schubert indexing $\sigma_{a_1,a_2}$, namely $\sigma_1(n-1)=\sigma_{1,1}$ and $\sigma_{k_i}(n)=\sigma_{k-1,0}$.
The sum of all the first indices $a_1$ appearing is
\[1\cdot(n-w(\pi))+\sum_{i=1}^{\ell(\pi)}(k_i-1) = n-w(\pi)+w(\pi)-\ell(\pi) = n-\ell(\pi) \le n-2.\]
The bound $c_1\le a_1+b_1$ in \eqref{computing_product} 
implies that when $\gamma(\pi)$ is expanded as a sum of Schubert classes,
none of the $\sigma_{d_1,d_2}$ arising in the computation have $d_1>n-2$, and therefore we never replace any symbol $\sigma_{d_1,d_2}$ by~0.
That is, $\gamma(\pi)$ may be computed in $\Lambda$ and we have a well-defined polynomial $\widehat s(\gamma(\pi))\in\mathbb Z[x,y]$.

We now conclude that $\widehat s\circ\gamma$ is injective, and therefore $\gamma$ is also.
This is because the factors $\widehat s(\sigma_{1,1})=xy$ and $\widehat s(\sigma_{k-1,0})=x^{k-1}+x^{k-2}y+\cdots+y^{k-1}$ are multiplicatively independent, with the exclusion of $\widehat s(\sigma_{0,0})=1$; 
factors of this last kind arise from parts $k_i=1$ in~$\pi$, and the number of parts $k_i=1$ is determined by $w(\pi)$ and the other parts.

To see the multiplicative independence, when we set $y=1$, the $k$\/th cyclotomic polynomial in~$x$ is a factor of $\widehat s(\sigma_{k-1,0})|_{y=1}$ but none of the $\widehat s(\sigma_{k'-1,0})|_{y=1}$ for $k'<k$,
and $\widehat s(\sigma_{1,1})|_{y=1}=x$ is not a factor of any $\widehat s(\sigma_{k-1,0})|_{y=1}$.
\fin
\end{dem}

Using the fact that $\sigma_{a_1,a_2} = (\sigma_{1,1})^{a_2}\cdot\sigma_{a_1-a_2,0}$,
the reader may check that the image of $\gamma$ is the set of all possible monomials in the Schubert classes $\{\sigma_{k}(w)\::\:1\leq k<w\leq n\}\subset H^{*}(\Gr (2,n),\mathbb{Z})$.

Since $\sigma_1(n)=\sigma_{0,0}=1\in H^{*}(\Gr (2,n),\mathbb{Z})$,
equation~\eqref{eq:def gamma} can be rewritten as in the next corollary.
\begin{cor}\label{cor_map}
If $\pi=(k_1,\ldots,k_t,1^m)$ with $k_i>1$ and $t\geq1$, then
\begin{align*}
\gamma\big((k_1,\ldots,k_t,1^m)\big)&=\sigma_1(n-1)^{n-w(\pi)}\prod_{k_i>1}\sigma_{k_i}(n)
\\&=\gamma\big((1^{n-1})\big)^{n-w(\pi)}\prod_{k_i>1}\gamma\big((k_i,1^{n-k_i})\big).
\end{align*}
\end{cor}

%

\begin{exa}
\label{Ex_pour_deux_part}
Let us pick $\pi=(k_1,k_2,1^{l-k_1-k_2})$, with $w(\pi)=l<n$. According to the previous Corollary, we need to solve $m_1+m_2=n+m$, with $m$ being $l-k_1-k_2$. One possibility is $m_1=n-k_1$ and $m_2=l-k_2$, so we have $\gamma(\pi)=\sigma_{k_1}(n)\,\sigma_{k_2}(l)=\gamma((k_1,1^{n-k_1}))\cdot \gamma((k_2,1^{l-k_2}))$.
\end{exa}

The next corollary resumes the notation of the end of Section~\ref{Section_T_Primitive}.

\begin{cor}
\label{cor_char_orb}
Let $\lambda\in H_*(\Gr(2,n),\mathbb{Z})$ be the class of a torus orbit. Then there exists a unique element $M\in M_n^2$ such that, for $y\in\Gr(2,n)^T$, we have $\theta(y)=\lambda$ if and only if $y\in\mathcal{G}_M(\mathbb{C})$. 
If we denote $\lambda$ by $\lambda(M)$, then $\mathcal{G}(2,n,\lambda(M))=\mathcal{G}_M(\mathbb{C})^{\sqcup\#(\mathfrak{S}_n\cdot M)}$.
\end{cor}
\begin{dem}
By Theorem \ref{teo_orbits}, we have $\mathcal{G}_M(\lambda(M))=\bigsqcup_N\mathcal{G}_N(\mathbb{C})$, where the union is indexed by the elements  $N\in\mathcal{M}_n^2$ such that $\lambda(N)=\lambda(M)$. By the injectivity in Theorem~\ref{thm_map} this happens if and only if $N\in\mathfrak{S}_n\cdot M$. So  ${M}_n^2(\lambda)=\{M\}$ (because the closure is the toric variety induced by the corresponding matroid polytope) and  by  Corollary \ref{cor_param_space} we have 
\begin{equation*}
    \mathcal{G}(2,n,\lambda)=\bigsqcup_{M\in {M}_n^2(\lambda)}\mathcal{G}_M(\lambda)^{\sqcup\#(\mathfrak{S}_n\cdot M)}=\mathcal{G}_M(\mathbb{C})^{\sqcup\#(\mathfrak{S}_n\cdot M)},
\end{equation*}
which ends the proof.\fin
\end{dem}

\subsection{Computing cohomology classes of orbits}
\label{Section_Computing}
This section is about the computation of $\mathbb{V}_{\pi}(*)$ where $*\in\mathcal{G}_\pi(\mathbb{C})$ is a closed point. 
Recall from Corollary \ref{cor_char_orb} that  if $\lambda\in H_*(\Gr(2,n),\mathbb{Z})$ is the class of a torus orbit, then there exists a unique element $\pi\in \Pi_n$ such that $\lambda=\lambda(\pi)=h(\mathbb{V}_{\pi}(*))$ for $*$ a closed point of~$\mathcal{G}_\pi(\mathbb{C})$.

If $\pi$ is rigid, then $\mathbb{V}_{\pi}(*)=\overline{G_\pi}$, so by Theorem~\ref{thm_map}
$[\mathbb{V}_{\pi}(*)]=\gamma(\pi)$, which is the Poincar\'e dual of $\lambda(\pi)$. 
We move on to the question of the  expression of $[\mb{V}_\pi(*)]$ when $\pi\in\Pi_n$ is not rigid, that is, $\ell(\pi)\geq4$. We will show that $[\mb{V}_\pi(*)]$ may be expressed as  $\sum_im_i[\overline{G}_{\pi_i}]$, where $m_i\in\mathbb{Z}_{>0}$ and each $\pi_i$ is rigid. 
We now prove this
by introducing an algorithm to compute effectively the non-rigid matroidal classes $[\mathbb{V}_{\pi}(*)]$. In particular, we can encode a regular matroidal subdivision $\{\Delta_{M_1},\ldots,\Delta_{M_k}\}$ of $\Delta_M$ in $\mathcal{M}^d_n$ as (the combinatorial type of) a  $(d-1)$-dimensional tropical linear space in the tropical projective space $\mathbb{TP}^{n-1}$. Thus the key for this procedure is Corollary \ref{cor_computation_orbits}, which says that we need to choose a combinatorial type of $1$-dimensional tropical linear space in $\mathbb{TP}^n$ corresponding to a maximal matroidal subdivision of $\Delta_\pi$.

\begin{defi}
 A \textbf{model} $L_\pi$ for $\pi=\prod_{i=1}^tk_i^{\alpha_i}$  is a pair $L_\pi=(T,m)$, where \begin{enumerate}
\item $T=(V,E)$ is a 3-valent tree with $\ell(\pi)$ vertices of degree one,
\item $m$ assigns, to every vertex $v$ of degree one, a part $m(v)=k_i$ of $\pi$.
\end{enumerate}
\end{defi}

Thus $T$ has $\ell(\pi)-2$ vertices of degree 3, and we consider  $T=(V,E)$ as the gluing of the $\ell(\pi)-2$ trivalent stars. We now turn each one of these stars into a cohomology class by assigning to  each one of these vertices $v$ of degree 3,  a partition $\pi(v)$ with $w(\pi(v))=w(\pi)$ and  $\ell(\pi)=3$ as follows. Let $C_1(v),C_2(v)$ and $C_3(v)$ be the 3 connected components of the graph $T\setminus\{v\}$. We denote by $\pi(v)$ the partition $m(C_1(v))m(C_2(v))m(C_3(v))$, where $m(C_i(v))=\sum_{u\in C_i(v)}m(u)$.

\begin{teo}
\label{teo_decompodecompo}
Let $\pi\in \Pi_n$ be non-rigid and let $L_\pi=(T,m)$ be a model for $\pi$. The orbit class $\lambda(\pi)$ is the Poincar\'e dual of  
\begin{equation}
\label{decompodecompo}
[\mb{V}_\pi(*)]=[\overline{T_\pi}]=\sum_{v\in L_\pi}\gamma(\pi(v)).
\end{equation}
The expression \eqref{decompodecompo} is independent of the choice of the model $L_\pi$.
\end{teo}
\begin{dem}
For $\pi$ nonrigid, we will see that a  model $L_\pi=(T,m)$ encodes a matroidal subdivision of 
$\Delta_\pi=\Delta(M_\pi)$ with  elements  $\{\Delta_{\pi(v)}\::\:v\in L_\pi\}$. We need to show that $\Delta_{\pi(v)}\cap \Delta_{\pi(w)}$ is a matroid for every $v\neq w$. After re-labelling, we can suppose that   $C_1(w),C_2(w)\subset C_3(v)$, $C_1(v),C_2(v)\subset C_3(w)$. We have that $(C_1(v)\cup C_2(v))\cap(C_1(w)\cup C_2(w))=\emptyset$.

Let $\pi=\{I_1,\ldots,I_{\ell(\pi)}\}$ be the partition of $[w(\pi)]$ induced by $\pi$ as in \eqref{particular_partition}. So, let $I(v)=\{i\in [w(\pi)]\::\:i\in I_k,\:k\in C_1(v)\cup C_2(v)\}$ and we define analogously $I(w)$. Then $\Delta_{\pi(v)}\cap \Delta_{\pi(w)}=U_{I(v)}^1\oplus U_{I(w)}^1.$

 It follows from Corollary \ref{cor_computation_orbits} that $$c_{w(\pi)-1}(M_\pi)=\sum_{v\in L_\pi}c_{w(\pi)-1}(M_{\pi(v)})=\sum_{v\in L_\pi}\gamma(\pi(v)),$$ 

and that this is independent from the choice of the model.\fin 
\end{dem}

\begin{exa}Note that $\pi=(k,1^{m-k})$ is non-rigid for $k=1,\ldots,m-3$ and $m\geq4$. Using the caterpillar model for $\pi$, we have:
\begin{equation*}
\lambda(\pi)=\sum_{i=k}^{m-2}[\overline{G}_{(i,m-i-1,1)}]
=\sum_{i=k}^{m-2}\delta(\sigma_{i}(n)\cdot\sigma_{m-i-1}(m)).
\end{equation*}
This last part comes from Example \ref{Ex_pour_deux_part}. The partition $(i,m-i-1,1)$ equals $\pi=(k_1,k_2,1^{l-k_1-k_2})$ with $k_1=i$, $k_2=m-i-1$. Then $l=m$, and from the expression $\gamma(\pi)=\sigma_{k_1}(n)\,\sigma_{k_2}(l)$, we get the result.
Finally, to compute each product $\sigma_{i}(n)\cdot\sigma_{m-i-1}(m)$, we use the change of coordinates $a_1=k+n-m-1$, $a_2=n-m$, and the formula \eqref{computing_product}.
\end{exa}

\section{Link with Euler-Chow series}
\label{Section_Link}
In this section we discuss how the foregoing results can be applied
to the computation of Euler-Chow series of Grassmannians of lines.
We complete this computation for $\mb{G}(2,4)$, which has been done in a different way in \cite{eli-lim}.
For $n>4$, the computation can be reduced to enumerating cycles within the spaces $\mathcal{G}_M$.
Carrying out these enumerations requires further technology beyond that in this paper and we leave it as an open problem.

We begin by defining Chow varieties and Euler-Chow series.
Let $X$ be a smooth projective complex variety and $p$ an integer.
Each effective algebraic $p$-cycle $Z=\sum_in_iY_i$ on~$X$ has a homology class $[Z]\in H_{2p}(X,\mathbb{Z})$, defined in terms of its summands as $[Z]=\sum_in_i[Y_i]$. 

The \textbf{Chow variety} $C_{p,\lambda}(X)$
parametrizes effective algebraic $p$-cycles in $X$ of a fixed homology class $\lambda\in H_{2p}(X,\mathbb{Z})$.
In the case where $H_{2p}(X,\mathbb{Z})\cong\mathbb{Z}$ has a unique generator $\lambda$ that is the class of an algebraic variety,
we abbreviate $C_{p,e}(X)=C_{p,e\lambda}(X)$.
We write 
\[C_{p,*}(X) := \bigsqcup_{\lambda\in H_{2p}(X,\mathbb{Z})} C_{p,\lambda}(X)\] 
for the set of all effective algebraic $p$-cycles in $X$.
We will only need to regard $C_{p,*}(X)$ as a set of points, rather than a scheme.
$C_{p,*}(X)$ is a monoid under sum of cycles.

For fixed integers $d,n,p$ with $0\leq p\leq 2n-4$, 
the $p$-dimensional \textbf{Euler-Chow series} of $\mb{G}(2,n)$ is defined as 
\begin{equation}\label{eq:ECS}
E_p(\mb{G}(d,n))=\sum_{\lambda\in H_{2p}(X,\mathbb{Z})}\chi(C_{p,\lambda}(\mb{G}(d,n)))\cdot t^\lambda\in\mathbb{Z}[\![M]\!],
\end{equation}
where 
$\chi$ is the topological Euler characteristic.

The $T$-action on $\mb{G}(d,n)$ induces a $T$-action on $C_{p,\lambda}(\mb{G}(d,n))$.
We showed in Theorem~\ref{V(Y)_2} and the discussion following that
every nonempty $C_{p,\lambda}(\mb{G}(d,n))$ has $T$-fixed points, 
namely effective algebraic $p$-cycles on $\mb{G}(d,n)$ of homology class $\lambda$ supported on $p$-dimensional subvarieties $Y\in \mb{G}(d,n)^T$.
By \cite{law&yau-hosy} we have
\[\chi(C_{p,\lambda}(\mb{G}(d,n)))=\chi(C_{p,\lambda}(\mb{G}(d,n))^T).\]
By decomposition into irreducible components we can write $C_{p,\lambda}(\mb{G}(d,n))^T$
in terms of the parameter spaces $\mathcal{G}(d,n,\lambda)$ of Corollary~\ref{cor_param_space}.
In this way, a sufficiently good understanding of subvarieties of the spaces $\mathcal{G}(d,n,\lambda)$
will enable the computation of $E_p(\mb{G}(d,n))$.
A key fact for the case $d=2$ is Theorem~\ref{thm_geometry},
identifying $\mathcal{G}(2,n,\lambda)$ as isomorphic to some $\mathcal{M}_{0,\ell}$.

Let us start with the case $\mb{G}(2,4)$.
For dimensions $p\ne3$,
the set of $T$-fixed subvarieties of dimension~$p$ in $\mb{G}(2,4)$ is finite:
there are 6 fixed points, 12 fixed curves, 11 fixed surfaces, 
and $\mb{G}(2,4)$ itself is the only fixed subvariety of dimension~4.
In these dimensions we simply have
\[E_p(\mb{G}(2,4))=\prod_{\stackbounds{Y\in \mb{G}(2,4)^T}{\dim Y=p}}\frac1{1-t^{[Y]}}.\]

So we focus on computing $E_3(\mb{G}(2,4))$, the 3-dimensional Euler--Chow series of $\mb{G}(2,4)$.
Since $H_6(\mb{G}(2,4),\mathbb{Z})\cong\mathbb{Z}\cdot[\Sigma_1(4)]$, this is $$E_3(\mb{G}(2,4))=\sum_{d\in\mathbb{Z}_{\geq1}}\chi(C_{3,d}(\mb{G}(2,4)))\cdot t^d,$$
where we recall the abbreviation $C_{3,d}(\mb{G}(2,4))=C_{3,d[\Sigma_1(4)]}(\mb{G}(2,4))$.

\begin{teo}\label{C_3d G24}
For all $d\geq1$, we have $\chi(C_{3,d}(\mb{G}(2,4)))=\binom{5+d}{5}-\binom{3+d}{5}=\frac1{12}(d+3)(d+2)^2(d+1)$.
\end{teo}

\begin{dem}
Following our strategy we compute $\chi(C_{3,d}(\mb{G}(2,4)))=\chi(C_{3,d}(\mb{G}(2,4))^T)$. 
All $T$-fixed 3-dimensional subvarieties are torus orbits, whose classes are $[\overline{T_{2\cdot1^2}}]=[\Sigma_1(4)]$ and $[\overline{T_{1^4}}]=2[\Sigma_1(4)]$.
There are $6$ $T$-invariant representatives $S_1,\ldots,S_6$ of the class $[\Sigma_1(4)]$, and an $\mathcal{M}_{0,4}$ worth of representatives $V_i$ of the class $2[\Sigma_1(4)]$.  
A cycle in $C_{3,d}(\mb{G}(2,4))$ has the form $n_1S_1+\cdots +n_6S_6+m_1V_1+\cdots +m_rV_r$,
where the $n_i,m_j\geq0$ satisfy $\sum_in_i+2\sum_jm_j=d$.

Now the part $m_1V_1+\cdots +m_sV_s$ of the cycle projects to an effective divisor of degree $m=m_1+\cdots+m_r$ on $\mathcal{M}_{0,4}$, and we know that the space of such divisors is isomorphic to $\mathcal{M}_{0,4}^{(m)}$, the $m$-th symmetric power of $\mathcal{M}_{0,4}$. So the Euler characteristic of the space of such divisors is $\chi(\mathcal{M}_{0,4}^{(m)})$. 
For a finite CW-complex $Z$, we have $\chi(Z^{(m)})=(-1)^m\binom{-\chi(Z)}{m}$. 
$\mathcal{M}_{0,4}$ is not a finite CW-complex but is homotopy-equivalent to a finite CW-complex $X$, the wedge of two circles,
inducing a homotopy equivalence between $\mathcal{M}_{0,4}^{(m)}$ and $X^{(m)}$.
It follows that $\chi(\mathcal{M}_{0,4}^{(m)})=(-1)^m\binom{1}{m}$, which is zero if $m>1$.

Thus we only need to take into account cycles of the form $n_1S_1+\cdots +n_6S_6$ subject to $\sum_in_i=d$, or $n_1S_1+\cdots +n_6S_6+V_1$ subject to $\sum_in_i=d-2$. The contribution of each of the $\binom{5+d}{5}$ cycles of the first form is 1, and the contribution of each of the $\binom{5+d-2}{5}$ cycles of the second form is $-1$.
\fin
\end{dem}

We close with remarks on the next case, $\mb{G}(2,5)$.
If we wanted to write down all Euler-Chow series for $\Gr(2,5)$ completely, 
the only part of the computation that does not have a counterpart in the previous section is the handling of (integral) $T$-invariant hypersurfaces $Y\subset \Gr(2,5)$ such that $Y=\overline{Y\cap G_{1^5}}$.
For such $Y$ there exists $d\in\mathbb{N}$ such that 
$[Y]=dH$, where $H=H_{\Gr(2,5)}\in H_{10}(\Gr(2,5),\mathbb{Z})$
is the class of a hyperplane section in the Pl\"ucker embedding, 
and the task would be to find the Euler characteristic of the space of~$Y$ for each fixed $d$.

Every such $Y$ is of the form $Y=\mathbb{V}_{1^5}(y)$ for some integral curve $y\in\mathcal{G}_{1^5}\cong\mathcal{M}_{0,5}$,
so the problem becomes one of curve counting in $\mathcal{M}_{0,5}$.
Using the presentation of the Deligne--Mumford compactification $\overline{\mathcal{M}_{0,5}}$ as isomorphic to the blowup $\text{Bl}_{\{p_1,\ldots,p_4\}}\mathbb{P}^2$,
this can be reduced to counting curves in $\mathbb{P}^2$
according to their degree and intersection multiplicity with the four blow-up centers.

\vspace{.3cm}
\par\noindent
\textsc{Instituto de Matem\'aticas, Universidad Nacional Aut\'onoma de M\'exico}
\par\noindent
e-mail address: \verb+javier@im.unam.mx+
\par\bigskip\noindent
\textsc{School of Mathematical Sciences, Queen Mary University of London
}
\par\noindent
e-mail address: \verb+a.fink@qmul.ac.uk+
\par\bigskip\noindent
\textsc{Centro de Investigaci\'on en Matem\'aticas, A.C., Jalisco S/N, Col. Valenciana CP: 36023 Guanajuato, Gto, M\'exico.
}
\par\noindent
e-mail address: \verb+cristhian.garay@cimat.mx+
\end{document}